\documentclass[12pt,reqno,a4wide]{article}
\usepackage{makeidx}
\usepackage{comment}
\usepackage{mathptmx}
\usepackage{a4wide}
\usepackage{amssymb}
\usepackage{amsfonts}
\usepackage{amsthm}
\usepackage{amsmath}
\usepackage{dsfont}
\usepackage{graphicx}
\usepackage{shadow}
\usepackage[all]{xy}
\usepackage{enumerate}
\usepackage{mathrsfs}
\usepackage{upgreek}
\usepackage{stmaryrd}
\usepackage[utf8]{inputenc}
\usepackage{array}
\usepackage{tikz-cd}
\usepackage{thumbpdf}
\usepackage[colorlinks=true,pagebackref]{hyperref}
\usepackage{textcomp}
\usepackage{chapterbib}
\usepackage{graphics}
\usepackage{longtable}
\usepackage{epsf}
\usepackage{bm}
\usepackage{adjustbox}

\newtheorem{theorem}{Theorem}[section]
\newtheorem{lem}[theorem]{Lemma}
\newtheorem{thm}[theorem]{Theorem}
\newtheorem{prop}[theorem]{Proposition}
\newtheorem{cor}[theorem]{Corollary}
\theoremstyle{definition}
\newtheorem*{Beweis}{Proof}
\newtheorem{defn}[theorem]{Definition}
\newtheorem{ex}[theorem]{Example}

\newtheorem{rem}[theorem]{Remark}
\newtheorem{rems}[theorem]{Remarks}
\newtheorem{punto}[theorem]{}

\swapnumbers
\CompileMatrices
\input{tcilatex}
\begin{document}

\title{Zariski-like Topologies for Lattices with Applications to Modules over Commutative Rings
\thanks{MSC2010: Primary 06A15; Secondary 16D10, 13C05, 13C13, 54B99. \newline
Key Words: Topological lattices; Prime Modules, First Submodules, Strongly
Hollow Submodules, Zariski Topology, Dual Zariski Topology, Zariski
Topology; Dual Zariski Topology}}
\author{ $%
\begin{array}{ccc}
\text{Jawad Abuhlail}\thanks{%
Corresponding Author; Email: abuhlail@kfupm.edu.sa. \newline
The authors would like to acknowledge the support provided by the Deanship
of Scientific Research (DSR) at King Fahd University of Petroleum $\&$
Minerals (KFUPM) for funding this work through projects No. RG1213-1 $\&$
RG1213-2} &  & \text{Hamza Hroub}\thanks{%
The paper is extracted from the Ph.D. dissertation of Dr. Hamza Hroub under
the supervision of Prof. Jawad Abuhlail} \\
\text{Department of Mathematics and Statistics} &  & \text{Department of
Mathematics } \\
\text{King Fahd University of Petroleum }\&\text{ Minerals} &  & \text{King
Saud University} \\
\text{31261 Dhahran, KSA} &  & \text{11451 Riyadh, KSA}%
\end{array}%
$}
\date{\today }
\maketitle

\begin{abstract}
We study Zariski-like topologies on a proper class $X\varsubsetneqq L$ of a
complete lattice $\mathcal{L}=(L,\wedge ,\vee ,0,1)$. We consider $X$ with
the so called \emph{classical Zariski topology} $(X,\tau ^{cl})$ and study
its topological properties (\emph{e.g.} the separation axioms, the
connectedness, the compactness) and provide sufficient conditions for it to
be \textit{spectral}. We say that $\mathcal{L}$ is $X$\emph{-top} iff%
\begin{equation*}
\tau :=\{X\backslash V(a)\mid a\in L\},\text{ where }V(a)=\{x\in L\mid a\leq
x\}
\end{equation*}%
is a topology. We study the interplay between the \textit{algebraic
properties} of an $X$-top complete lattice $\mathcal{L}$ and the \textit{%
topological properties} of $(X,\tau ^{cl})=(X,\tau ).$ Our results are
applied to several spectra which are proper classes of $\mathcal{L}%
:=LAT(_{R}M)$ where $M$ is a left module over an arbitrary associative ring $%
R$ (e.g. the spectra of \emph{prime}, \emph{coprime}, \emph{fully prime}
submodules) of $M$ as well as to several spectra of the dual complete
lattice $\mathcal{L}^{0}$ (\emph{e.g.} the spectra of \emph{first}, \emph{%
second} and \emph{fully coprime} submodules of $M$).
\end{abstract}

\section*{Introduction}

The spectrum $Spec(R)$ of prime ideals of a commutative ring $R$ attains the
so called \emph{Zariski topology} in which the closed sets are the varieties
\begin{equation*}
\{V(I)\mid I\in Ideal(R)\},\text{ where }V(I)=\{P\in Spec(R)\mid I\subseteq
P\}.
\end{equation*}%
This topology is compact, $T_{0}$ but almost never $T_{2},$ and the closed
points correspond to the maximal ideals. The Zariski topology proved to be
very important in two main aspects: in Algebraic Geometry and in Commutative
Algebra. In particular, it provided an efficient tool for studying the
algebraic properties of a commutative ring $R$ by investigating the
corresponding topological properties of $Spec(R)$ \cite{AM1969}.

Motivated by this, there were many attempts to define Zariski-like
topologies on the spectra of \emph{prime-like }submodules of a given left
module $M$ over a (not necessarily commutative) ring $R.$ This resulted at
the first place in several different notions of prime submodules of $_{R}M$
which reduced to the notion of a prime ideal for the special case $M=$ $R,$
a commutative ring (\emph{e.g.} \cite{Wij2006}). The work in this direction
was almost limited to studying these prime-like submodules and their duals
(the coprime-like submodules) as well as to the families of prime ideals
corresponding to them from a purely algebraic point of view. One of the
obstacles was that not every module $M$ over a (commutative) ring $R$ has
the property that $Spec(M)$ attains a Zariski-like topology: the proposed
\emph{closed varieties} $\{V(N)\mid N\in LAT(_{R}M)\}$ are not necessarily
closed under finite unions. Modules for which this last condition is
satisfied were investigated, among others, by R. L. McCasland and P. F.
Smith (e.g. \cite{MMS1997}, \cite{M}) and called \emph{top modules}.
However, even such modules were studied from a purely algebraic point of
view and the associated Zariski-like topologies were not well studied till
about a decade ago. In \cite{Abu2006}, Abuhlail introduced a Zariski-like
topology on the spectrum of \emph{fully coprime subcomodules} of a given
comodule $M$ of a coring $\mathcal{C}$ over an associative ring $R$ and
studied the interplay between the algebraic properties of $M$ and the
topological properties of that Zariski-like topology (see also \cite{Abu2008}).

Later, in a series of papers (\cite{Abu2011-a}, \cite{Abu2011-b}, \cite{Abu}),
Abuhlail introduced and investigated several Zariski-like topologies for
a module $M$ over an arbitrary associative ring $R.$ These investigations
showed that all the (co)prime spectra considered fall in two main classes
with several common properties for the spectra in each class. Moreover,
these two classes were \emph{dual }to each other in some sense. This led
Abuhlail and Lomp (\cite{AbuC0}, \cite{AbuC}) to investigate such topologies
for a general complete lattice $\mathcal{L}:=(L,\wedge ,\vee ,0,1)$ and a
proper subset $X\subseteq L\backslash \{1\}.$ Their main work was on
characterizing the so called $X$\emph{-top lattices} (i.e. $\mathcal{L}$ for
which the closed varieties $V(a):=\{x\in X\mid a\leq x\}$ are closed under
finite unions). In addition to the fact that this approach provided a
general framework, it allowed obtaining results on the dual lattice $%
\mathcal{L}^{0}:=(L,\wedge ^{0},\vee ^{0},0^{0},1^{0})=(L,\vee ,\wedge ,1,0)$
and $X\subseteq L\backslash \{0\}$ for free.

This paper consists of two sections. After providing some basic
definitions and preliminaries, we study in Section 1 Zariski-like
topologies for complete lattices using a different approach. Fix a complete
lattice $\mathcal{L}=(L,\vee ,\wedge ,1,0),$ a subset $X\subseteq
L\backslash \{1\}$ and $\tau :=\{X\backslash V(a)\mid a\in L\}.$ Inspired by
the work of Behboodi and Haddadi \cite{BH2008-a}, \cite{BH2008-b} on the
lattice $LAT(_{R}M)$ of submodules of a given module $M$ over a ring $R,$
and instead of restricting our attention to $X$\emph{-top lattices} (i.e.
those for which $(X,\tau )$ is a topology), we consider $X$ with the \emph{%
classical Zariski topology} $(X,\tau ^{cl})$ which is constructed on $X$ by
considering $\tau $ as a \emph{subbase} and the \emph{finer patch topology} $%
(X,\tau ^{\mathfrak{fp}})$ which has a subbase $\mathcal{B}:=\{V(a)\cap
X\backslash V(b)\mid a,b\in L\}.$ Indeed, $(X,\tau ^{cl})\leq (X,\tau ^{%
\mathfrak{fp}})$ and $(X,\tau )=(X,\tau ^{cl})$ if and only if $\mathcal{L}$
is $X$-top. In the special case when $\mathcal{L}$ is an $X$-top lattice, we
not only apply the results obtained on $(X,\tau ^{cl})$, but obtain also
other interesting results especially on the interplay between the algebraic
properties of $\mathcal{L}$ and the topological properties of $(X,\tau ).$

In Proposition \ref{Proposition 3.11}, we prove a stronger version of the
converse of \cite[Proposition 2.7]{AbuC} and conclude in Corollary \ref%
{Corollary 3.12} that in case $\mathcal{L}$ is an $X$-top lattice: $%
A\subseteq X$ is irreducible if and only if $I(A):=\bigwedge\limits_{x\in
A}x $ is (strongly) irreducible in the sublattice $(\mathcal{C}(L),\wedge )$
of radical elements of $\mathcal{L}.$ This fact was the key in the proofs of
several results including Theorem \ref{Proposition 3.48}.
It is worth mentioning that Theorem \ref{Proposition 3.48}
recovers several results of Abuhlail on such 1-1 correspondences for $%
\mathcal{L}=LAT(_{R}M)$ (\emph{e.g.} \cite{Abu}, \cite{Abu2011-a}, \cite%
{Abu2011-b}) and Abuhlail/Lomp \cite{AbuC} as special cases (some of these
results are recovered under conditions weaker than those assumed in the
original results for the different spectra of modules).

In Theorem \ref{Theorem 3.16 000}, we prove that the class $Max(X)$ of
maximal elements of $X$ coincides with the class of $Max(\mathcal{C}(%
\mathcal{L}))$ of maximal radical elements. This yields, assuming that $%
\mathcal{C}(\mathcal{L})$ satisfies the so called \emph{complete max property%
}, that $(X,\tau ^{cl})$ is discrete if and only if $(X,\tau ^{cl})$ is $%
T_{1}.$ This result generalizes \cite[Theorem 5.34]{Abu}, \cite[Theorem 4.28]%
{Abu} and \cite[Theorem 3.46]{Abu2011-a}.

A topological space $T$ is said to be \emph{spectral} \cite{H1969} iff $T$
is homeomorphic to $Spec(R),$ the prime spectrum of a commutative ring $R,$
with the Zariski topology. Hochster \cite{H1969} characterized such spaces by giving
sufficient and necessary conditions on a topological space to be spectral.
We observe in Proposition \ref{Proposition 3.26} that if the finer patch topology
$(X,\tau ^{\mathfrak{fp}})$ is compact, then the classical Zariski topology $(X,\tau ^{cl})$ is spectral.
Sufficient conditions for $(X,\tau ^{\mathfrak{fp}})$ to be compact were provided in
Theorems \ref{Theorem 3.33} and \ref{Theorem 3.36}. Example \ref{Example
3.42} provides several spectra of modules which are shown to be spectral by
Theorem \ref{Theorem 3.33}.

In Section 2, we restrict our investigations to $X$-top lattices $\mathcal{L}%
=(L,\vee ,\wedge ,1,0)$ where $X\subseteq L\backslash \{1\}.$ We investigate
the interplay between the algebraic properties of $\mathcal{L}$ and the
topological space $(X,\tau )=(X,\tau ^{cl}).$ Several types of compactness
and connectedness of $(X,\tau )$ are studied in Theorem \ref{Theorem 3.45}.
For examples of such an interplay.

The results in Section 1 are applied to the complete lattice $LAT(_{R}M):=(%
\mathcal{L}(M),\cap ,+,0,M)$ of submodules of a left module $M$ over an
associative ring $R.$ In a series of examples \ref{ex-p} - \ref{ex-f}, we
apply Theorem \ref{Proposition 3.48} to a number of spectra $X\subseteq
\mathcal{L}(M)\backslash M$ (or $X\subseteq \mathcal{L}(M)\backslash \{0\}$).

\section{Zariski-like Topologies for Lattices}

\bigskip

\subsection*{Lattices}

\begin{punto}
(\cite{G}) A \emph{lattice}
\index{lattice} $\mathcal{L}$ is a \emph{poset} $(L,\leq )$ closed under two
binary commutative associative and idempotent operations $\wedge $ (\emph{%
meet}) and $\vee $ (\emph{join}), and we write $\mathcal{L}=(L,\wedge ,\vee
) $. We say that a lattice $(L,\wedge ,\vee )$ is a \emph{complete lattice} iff
$\bigwedge\limits_{x\in H}x$ and $\bigvee\limits_{x\in H}x$ exist for any $%
H\subseteq L$. For two lattices $\mathcal{L}=(L,\wedge ,\vee )$ and $%
\mathcal{L}^{\prime }=(L^{\prime },\wedge ^{\prime },\vee ^{\prime }),$ a
\emph{homomorphism of lattices}
\index{homomorphism of lattices} from $\mathcal{L}$ to $\mathcal{L}^{\prime
} $ is a map $\varphi :L\longrightarrow L^{\prime }$ that preserves finite
meets and finite joins, \emph{i.e.}%
\begin{equation*}
\varphi (x\wedge y)=\varphi (x)\wedge ^{\prime }\varphi (y)%
\text{ and }\varphi (x\vee y)=\varphi (x)\vee ^{\prime }\varphi (y)\text{ }%
\forall x,y\in L.
\end{equation*}%
If $\mathcal{L}=(L,\wedge ,\vee ,0,1)$ and $\mathcal{L}^{\prime }=(L^{\prime
},\wedge ^{\prime },\vee ^{\prime },0^{\prime },1^{\prime })$ are complete
lattices, then a \emph{morphism of complete lattices }from $\mathcal{L}$ to $%
\mathcal{L}^{\prime }$ is a map $\varphi :L\longrightarrow L^{\prime }$ that
preserves arbitrary meets and arbitrary joins.
\end{punto}

\begin{punto}
\index{bounded lattice}%
\index{complete lattice} Let $\mathcal{L}=(L,\wedge ,\vee )$ be a lattice.
If $\mathcal{L}$ has a \emph{maximum element} $1$ and a \emph{minimum
element }$0,$ then $\mathcal{L}$ is called a \emph{bounded lattice} and we
write $\mathcal{L}=(L,\wedge ,\vee ,0,1).$ An element $x\in L\backslash
\{1\} $ is called \emph{maximal} in $\mathcal{L}$ iff $y=x$ or $y=1$
whenever $x\leq y;$ dually, an element $x\in L\backslash \{0\}$ is called
\emph{minimal} iff $y=x$ or $y=0$ whenever $y\leq x.$ Notice that every
complete lattice is bounded. We make the convention that $%
\bigwedge\limits_{x\in \emptyset }x=1$ and $\bigvee\limits_{x\in \emptyset
}x=0.$
\end{punto}

\begin{punto}
\index{dual lattice} For every lattice $\mathcal{L}=(L,\wedge ,\vee ),$
there is associated the \emph{dual lattice} $\mathcal{L}^{0}=(L,\wedge
^{0},\vee ^{0})$ where $\wedge ^{0}=\vee $ and $\vee ^{0}=\wedge .$ Indeed,
if $\mathcal{L}=(L,\wedge ,\vee )$ is a complete lattice, then the dual
lattice $\mathcal{L}^{0}$ is complete. Moreover, if $\mathcal{L}=(L,\wedge
,\vee ,0,1)$ is a bounded lattice, then the dual lattice $\mathcal{L}%
^{0}=(L,\wedge ^{0},\vee ^{0},0^{0},1^{0})$ is bounded with $0^{0}=1$ and $%
1^{0}=0.$
\end{punto}

\begin{ex}
\label{ideal}%
\index{$LAT(M)$}%
\index{$Ideal(R)$} Let $R$ be a ring.

\begin{enumerate}
\item $S=(Ideal(R),\cap ,+,R,0),$ where $Ideal(R)$ is the set of all
(two-sided) ideals of $R$ is a complete lattice.

\item For any left $R$-module $M,$ the set $LAT(M)=(\mathcal{L}(M),\cap
,+,M,0)$ is a complete lattice where $\mathcal{L}(M)$ is the class of all $R$%
-submodules of $M$.
\end{enumerate}
\end{ex}

\begin{punto}
\index{irreducible element}
\index{strongly irreducible element}
\index{hollow element}
\index{strongly hollow element} Let $\mathcal{L}=(L,\wedge ,\vee ,0,1)$ be a
complete lattice.

\begin{enumerate}
\item An element $x\in L\backslash \{1\}$ is said to be:

\emph{irreducible} \cite{AbuC0} iff for any $a,b\in L$ with $a\wedge b=x,$
we have $a=x$ or $b=x$;

\emph{strongly irreducible} \cite{AbuC0} iff for any $a,b\in L$ with $%
a\wedge b\leq x,$ we have $a\leq x$ or $b\leq x$.

We denote the set of strongly irreducible elements in $L$ by $SI(\mathcal{L}%
) $.

\item An element $x\in L\backslash \{0\}$ is said to be:

\emph{hollow }iff whenever for any $a,b\in L$ with $x=a\vee b,$ we have $x=a$
or $x=b;$

\emph{strongly hollow} \cite{AbuC0} iff for any $a,b\in L$ with $x\leq a\vee
b,$ we have $x\leq a$ or $x\leq b$.

We denote the set of strongly hollow elements in $L$ by $SH(\mathcal{L})$.

\item We say that $\mathcal{L}$ is

a \emph{hollow lattice}%
\index{hollow lattice} iff $1$ is hollow (\emph{i.e.} for any two elements $%
x,y\in L\backslash \{1\}$ we have $x\vee y\neq 1$);

a \emph{uniform lattice}%
\index{uniform lattice} iff $0$ is uniform (\emph{i.e.} for any two elements
$x,y\in L\backslash \{0\}$ we have $x\wedge y\neq 0$).
\end{enumerate}
\end{punto}

\subsection{$X$-top Lattices%
\index{$X-top$}}

From now on, we assume that $\mathcal{L}=(L,\wedge ,\vee ,0,1)$ is a
complete lattice.

\begin{punto}
Let $X\subseteq L\backslash \{1\}.$ For $a\in L,$ we define \textit{the
variety of }$a$%
\index{variety} as $V(a):=\{p\in X\mid a\leq p\}$ and set $V(\mathcal{L}%
):=\{V(a)\mid a\in L\}.$ Indeed, $V(\mathcal{L})$ is closed under arbitrary
intersections (in fact, $\bigcap_{a\in A}V(a)=V(\bigvee_{a\in A}(a))$ for
any $A\subseteq L$). The lattice $\mathcal{L}$ is called $X$\emph{-top}
(or a \textit{topological lattice} \cite{AbuC}) iff $V(\mathcal{L})$ is closed under finite
unions. The lattice $\mathcal{L}$ is called \emph{strongly} $X$\emph{-top}
iff $X\subseteq SI(\mathcal{L})$ \cite{AbuC}.%
\index{strongly $X-top$ lattice}
\end{punto}

\begin{punto}
\index{radical element} Let $X\subseteq L\backslash \{1\}.$ For any $%
Y\subseteq X,$ we set $I(Y):=\bigwedge\limits_{p\in Y}p$ and $%
\sqrt{a}:=I(V(a)).$ We say that $a$ is an $X$-\emph{radical element }iff $%
\sqrt{a}=a.$ The set of $X$-radical elements of $L$ is%
\index{$\mathcal{C}(\mathcal{L})$}%
\begin{equation*}
\mathcal{C}^{X}(L):=\{a\in L\mid
\sqrt{a}=a\}.
\end{equation*}%
When $X$ is clear from the context, we drop it from the above notation.
Notice that $\mathcal{C}(\mathcal{L})=(\mathcal{C}(L),\wedge ,\tilde{\vee},%
\sqrt{0},1)$ is a complete lattice, where $\tilde{\bigvee }Y:=IV(\bigvee
(Y)) $ for any $Y\subseteq \mathcal{C}(L),$ i.e. $\tilde{\bigvee }_{x\in Y}x=%
\sqrt{\bigvee\limits_{x\in Y}x}.$ It was proved in \cite[Theorem 2.2]{AbuC}
that $\mathcal{L}$ is an $X$-top lattice if and only if the map
\begin{equation*}
V:(\mathcal{C}(L),\wedge ,\tilde{\vee},1,\sqrt{0})\longrightarrow (\mathcal{P%
}(X),\cap ,\cup ,X,\emptyset ),\text{ }a\mapsto V(a)
\end{equation*}%
is an anti-homomorphism of lattices, that is
\begin{equation*}
V(a\wedge b)=V(a)\cup V(b)\text{ and }V(a\vee b)=V(a)\cap V(b)\text{ for all
}a,b\in \mathcal{C}(L).
\end{equation*}
\end{punto}

The following lemma appeared in \cite{AbuC} except for (2) which is clear.

\begin{lem}
\label{lemma 1.19}Let $X\subseteq L\backslash \{1\}.$ For any $x,y\in L$ and
$A,B\subseteq L$ we have:

\begin{enumerate}
\item $A\subseteq B \Rightarrow I(B)\leq I(A)$.

\item $V(x)\subseteq V(y) \Leftrightarrow \sqrt{y}\leq \sqrt{x}$. It follows
that $V(x) = V(y) \Leftrightarrow \sqrt{y} = \sqrt{x}$.

\item $V(x)=V(\sqrt{x})$.

\item $\bigcap\limits_{x\in A}V(x)=V(\bigvee\limits_{x\in A}(x))$.

\item $I\circ V\circ I = I$.

\item $V\circ I\circ V = V$.

\item $\mathcal{L}$ is $X$-top $\Longleftrightarrow V(x)\cup V(y)=V(x\wedge
y)$ for any $x,y\in \mathcal{C}(L).$
\end{enumerate}
\end{lem}

\begin{punto}
\label{tcl}Let $X\subseteq L\backslash \{1\}$ and set $\tau :=\{X\backslash
V(a)\mid a\in L\}.$ We define $\tau ^{cl}$%
\index{$\tau ^{cl}$} to be the topology constructed on $X$ by taking $\tau $
as a \emph{subbase}, that is $\tau ^{cl}$ is the set of all arbitrary unions
of finite intersections of elements in $\tau $, and is called the \emph{%
classical Zariski topology} on $X$%
\index{classical Zariski topology}. Moreover, $\mathcal{L}$ is $X$-top (%
\emph{i.e.} $\tau $ is closed under finite intersections) if and only if $%
\tau ^{cl}=\tau $.
\end{punto}

\begin{punto}
Let $Y\subseteq L\backslash \{0\}$. For any $a\in L$, we define \textit{the }%
\emph{dual variety}
\index{dual variety} $V^{0}(a):=\{q\in Y\mid q\leq a\}$ and set $V^{0}(%
\mathcal{L})=\{V^{0}(a)\mid a\in L\}$. We say that $\mathcal{L}$ is
\index{dual $Y-top$} \emph{dual }$Y$\emph{-top} iff the dual lattice $%
\mathcal{L}^{0}$ is a $Y$-top lattice. For any subset $A\subseteq Y$, we set
$H(A):=\bigvee_{q\in A}q$; also we set $%
\sqrt{a}^{0}:=H(V^{0}(a))$, and $\mathcal{H}(\mathcal{L})=\mathcal{C}^{Y}(%
\mathcal{L}^{0})$. \textit{The }\emph{dual classical Zariski topology}%
\index{dual classical zariski topology} $\tau ^{dcl}$ on $Y$ is constructed
by taking $\tau ^{0}:=\{Y\backslash V^{0}(a)\mid a\in L\}$ as a subbase for
this topology. With this process, one can dualize the results obtained in
this section for the (classical) Zariski-topology to results on the dual
(classical) Zariski topology.
\end{punto}

The following lemma recovers \cite[5.14 and 4.10]{Abu}, \cite[3.23]%
{Abu2011-a} and \cite[3.21]{Abu2011-b}.

\begin{lem}
\label{Lemma 3.3}Let $X\subseteq L\backslash \{1\}$ and assume that $%
\mathcal{L}$ is an $X$-top lattice. The closure of any $Y\subseteq X$ is
given by $%
\overline{Y}=V(I(Y))$.
\end{lem}

\begin{Beweis}
Let $Y\subseteq X.$ Notice that $\overline{Y}=V(a)$ for some $a\in L,$
whence $a\leq \bigwedge\limits_{p\in Y}p=I(Y)$ and so $V(I(Y))\subseteq V(a)=%
\overline{Y}$. On the other hand, $Y\subseteq V(I(Y))$ and so $\overline{Y}%
\subseteq V(I(Y))$.
\end{Beweis}

$\blacksquare$

\begin{punto}
A non-empty topological space $(T,\tau )$ is said
to be:

\begin{enumerate}
\item \emph{connected} iff $T$ is \emph{not} the union of two disjoint
non-empty open subsets (equivalently, $T$ is \emph{not} the union of two
disjoint non-empty closed sets).

\item \emph{hyperconnected} (or \emph{irreducible }\cite{B1966}) iff no two
non-empty open sets in $T$ are disjoint (equivalently, $T$ is \emph{not} the
union of two closed subsets).

\item \emph{ultraconnected} \cite{B1966} iff no two non-empty closed sets in
$T$ are disjoint.
\end{enumerate}
\end{punto}

\begin{punto}
Let $(T,\tau )$ be a topological space. A subset $%
A\subseteq T$ is called \emph{hyperconnected} \cite{B1966} (or \emph{%
irreducible}) iff $A$ is so when considered as a topological space w.r.t.
the relative topology induced from $(T,\tau )$ (equivalently, $A$ is
non-empty and for any two closed subsets $F_{1},$ $F_{2}$ in $T$ with $%
A\subseteq F_{1}\cup F_{2},$ we have $A\subseteq F_{1}$ or $A\subseteq F_{2}$%
). The empty set is \emph{not} considered to be irreducible. A closed subset
$F\subseteq T$ is said to have a \emph{generic point} $g\in T$ \cite{B1966}
iff $%
\overline{\{g\}}=F.$ The topological space $(T,\tau )$ is called \emph{sober}
iff every closed irreducible subset of $T$ has a \emph{unique} generic point.
\end{punto}

\begin{punto}
\index{irreducible component} A subset $A\subseteq T$ is irreducible if and
only if the closure $%
\overline{A}$ is irreducible. An \emph{irreducible component} \cite{B1966}
is an irreducible subset of $X$ which is not a proper subset of any
irreducible subset of $T$ (hence an irreducible component of $T$ is indeed a
closed subset).
\end{punto}

The following result generalizes \cite[3.2 and 3.3]{BH2008-a}.

\begin{prop}
\label{Proposition 3.4}Let $X\subseteq L\backslash \{1\}$ and consider $%
(X,\tau ^{cl})$.

\begin{enumerate}
\item For each $p\in X$, we have $\overline{\{p\}}=V(p)$.

\item $V(p)$ is irreducible $\forall p\in X$.

\item If $Y\subseteq X$ is closed, then $Y=\bigcup\limits_{p\in Y}V(p).$

\item If $\mathcal{L}$ is $X$-top, then for any closed subset $Y\subseteq X$
we have $Y=\bigcup\limits_{p\in Y}V(p)=V(\bigwedge\limits_{p\in Y}{p})$.
\end{enumerate}
\end{prop}

\begin{Beweis}
Consider $(X,\tau ^{cl})$.

\begin{enumerate}
\item Observe that $V(p)$ is closed in $(X,\tau ^{cl}),$ whence $\overline{%
\{p\}}\subseteq V(p)$. On the other hand, suppose that $\overline{\{p\}}%
=\bigcap_{i\in I}(\bigcup_{j=1}^{j=n_{i}}(V(x_{ij}))),$ where $x_{ij}\in L$.
Since $p\in \bigcup_{j=1}^{j=n_{i}}(V(x_{ij}))\hspace{0.3cm}\forall i\in I$,
it follows that $V(p)\subseteq \bigcup_{j=1}^{j=n_{i}}(V(x_{ij}))\hspace{%
0.3cm}\forall i\in I$. Therefore, $V(p)\subseteq \overline{\{p\}}$. Clearly,
$\overline{\{p\}}\subseteq $ $V(p),$ whence $\overline{\{p\}}=V(p)$.

\item Notice that $\{p\}$ is irreducible, whence $V(p)=\overline{\{p\}}$ is
irreducible.

\item Clear.

\item Let $Y\subseteq X$ be closed. It follows from (3) that $%
Y=\bigcup\limits_{p\in Y}V(p)\subseteq V(\bigwedge\limits_{p\in Y}{p}).$
Since $\mathcal{L}$ is assumed to be $X$-top, $Y=V(x)$ for some $x\in L$ and
so $x\leq \bigwedge\limits_{p\in Y}{p}$, whence $V(\bigwedge_{p\in Y}{p}%
)\subseteq V(x)=Y$. Consequently, $Y=\bigcup\limits_{p\in
Y}V(p)=V(\bigwedge\limits_{p\in Y}{p})$.
\end{enumerate}
\end{Beweis}

$\blacksquare$

\begin{ex}
\label{Example 3.5}Consider the complete ideal lattice $\mathcal{L}=(Ideal(%
\mathbb{Z}),\cap ,+,\mathbb{Z},0).$ Consider $X=Spec^{p}(\mathbb{Z})$, the
prime spectrum of $\mathbb{Z}.$ It is clear that $(X,\tau )$ is a
topological space (the usual Zariski topology on the spectrum of the
commutative ring $\mathbb{Z}$). Notice that for $Y:=Spec^{p}(\mathbb{Z}%
)\backslash \{0\},$ we have $\overline{Y}=V(I(Y))=V(\bigcap\limits_{P\in
Y}P)=V(0)=X\neq \bigcup\limits_{p\in Y}V(p)$. This example shows that \cite[%
Proposition 3.1]{BH2008-a} fails to hold even for domains. However, the
proof of Proposition \ref{Proposition 3.4} provides a correct proof \cite[%
Corollary 2.3]{BH2008-a} without using \cite[Proposition 3.1]{BH2008-a}.
\end{ex}

The following result recovers \cite[Proposition 3.8]{BH2008-a}, \cite[%
Proposition 3.24 (1)]{Abu2011-a} and \cite[Proposition 5.15 (i)]{Abu}.

\begin{prop}
\label{Proposition 3.6}Let $X\subseteq L\backslash \{1\}$ and consider $%
(X,\tau ^{cl})$.

\begin{enumerate}
\item $X$ is $T_{0}$.

\item Every finite closed irreducible subset of $X$ has a unique generic
point. If $X$ is finite, then $X$ is sober.
\end{enumerate}
\end{prop}

\begin{Beweis}
\begin{enumerate}
\item Let $p_{1},p_{2}\in X$ be such that $\overline{\{p_{1}\}} = \overline{%
\{p_{2}\}},$ whence $V(p_{1}) = V(p_{2})$ and it follows that $p_{1} = p_{2}$%
, which proves that $X$ is $T_0$ (notice that, in general, $(X, \tau)$ is $%
T_0$ if and only if $\overline{\{p_{1}\}} = \overline{\{p_{2}\}}
\Leftrightarrow p_{1} = p_{2}$).

\item In general, If $(X, \tau)$ is $T_0$, then every finite irreducible
closed subset has a unique generic point. To see this, suppose that $F$ is a
closed irreducible finite set that has no generic point. Pick $x_1\in F$,
whence $\overline{\{x_1\}}\neq F$ and so there is $x_2\in F\backslash
\overline{\{x_1\}}$. Observe that $\overline{\{x_1\}}\cup \overline{\{x_2\}}
\neq F$ as $F$ is irreducible. So, there is $x_3\in F\backslash (\overline{%
\{x_1\}}\cup \overline{\{x_2\}})$. by continuing this process, we
conclude that $F$ is infinite, a contradiction. The uniqueness of the
generic point follows directly from the fact that $T_0$.
\end{enumerate}
\end{Beweis}

$\blacksquare$

The following observation generalizes \cite[Proposition 2.3]{BH2008-a}.

\begin{rem}
\label{Remark 3.7}Let $X\subseteq L\backslash \{1\}.$ The following are
equivalent for $(X,\tau ^{cl}):$

\begin{enumerate}
\item $\mathcal{L}=\mathcal{C}(\mathcal{L})$.

\item For all $x_{1},x_{2}\in L$ with $V(x_{1})=V(x_{2}),$ we have $%
x_{1}=x_{2}.$
\end{enumerate}
\end{rem}

\begin{Beweis}
(1 $\Rightarrow $ 2) Suppose $V(x_{1})=V(x_{2})$ for some $x_{1},x_{2}\in L.$
It follows that $x_{1}\leq p,\forall p\in V(x_{2})$ whence $x_{1}\leq \sqrt{%
x_{2}}=x_{2}$. Similarly, $x_{2}\leq x_{1}$.

(2 $\Rightarrow $ 1) $\forall x\in L$ we have $V(x)=V(\sqrt{x}),$ whence $x=%
\sqrt{x}.$
\end{Beweis}

$\blacksquare$

\begin{punto}
\index{atomic}
\index{coatomic} \label{min}Let $X\subseteq L\backslash \{1\}$ and denote by
$Min(X)$ the set of minimal elements of $X$ and by $Max(X)$ the set of
maximal elements of $X.$ We say that $X$ is

\emph{atomic} iff for every $p\in X$ there is $q\in Min(X)$ such that $q\leq
p;$

\emph{coatomic} iff for every element $p\in X$ there is $q\in Max(X)$ such
that $p\leq q$.
\end{punto}

\begin{rems}
\label{Remark 3.9}Let $X\subseteq L\backslash \{1\}$ and consider $(X,\tau
^{cl})$.

\begin{enumerate}
\item If $X$ satisfies the DCC, then $X$ is atomic.

\item If $X$ is atomic, then there is a subset $A\subseteq X$ such that $%
X=\bigcup_{p\in A}V(p)$ with $V(p)$ and $V(q)$ are not comparable for any $%
p\neq q\in A$ (e.g. take $A=Min(X)$).

\item Let $X$ be atomic and $Min(X)$ finite. Then $X$ is irreducible if and
only if $Min(X)$ is a singleton. To see this, observe that $%
X=\bigcup\limits_{p\in Min(X)}V(p)$ with $p\nleq q$ for any $p\neq q$ are in
$Min(X)$. Clearly, $X$ is irreducible if and only if $Min(X)$ is a singleton.
\end{enumerate}
\end{rems}

\begin{rems}
\label{Remark 3.10}Let $X\subseteq L\backslash \{1\}$ with $0\in X$ and
consider $(X,\tau ^{cl}).$

\begin{enumerate}
\item If $F\subseteq X$ is closed and $0\in F$, then $F=X$. To prove this,
observe that $X=V(0)=%
\overline{\{ 0\} }\subseteq F$.

\item Every non-empty open subset of $X$ contains $0$. To see this, let $%
O\subseteq X$ be open. If $0\notin O,$ then $0\in F:=X\backslash O$. By (1),
$X\backslash O=X,$ i.e. $O=\emptyset $.

\item $X$ is irreducible since $Min(X)=\{0\}$, a singleton (see Remark \ref%
{Remark 3.9} (3)).
\end{enumerate}
\end{rems}

\noindent It was proved in \cite[Proposition 2.7]{AbuC}, that if $\mathcal{L}
$ is an $X$-top lattice and $A\subseteq X$ is such that $I(A)$ is
irreducible in $(\mathcal{C}(L),\wedge ),$ then $A$ is irreducible in $%
(X,\tau )$. The following result proves a stronger version of the converse.

\begin{prop}
\label{Proposition 3.11}Let $X\subseteq L\backslash \{1\}$ and assume that $%
\mathcal{L}$ is an $X$-top lattice. If $A\subseteq X$ is irreducible, then $%
I(A)$ is strongly irreducible in $(\mathcal{C}(L),\wedge )$.
\end{prop}

\begin{Beweis}
Suppose that $a\wedge b\leq I(A)$ for some $a,b\in \mathcal{C}(L)$. Now, $%
\overline{A}=V(I(A))\subseteq V(a\wedge b)\overset{\text{\cite[Theorem 2.2]%
{AbuC}}}{=}V(a)\cup V(b).$ Since $A$ is irreducible, $\overline{A}$ is also
irreducible, whence $\overline{A}\subseteq V(a)$ or $\overline{A}\subseteq
V(b)$. So, $a=I(V(a))\leq I(\overline{A})=I(V(I(A)))=I(A)$ or $b=I(V(b))\leq
I(\overline{A})=I(V(I(A)))=I(A)$.
\end{Beweis}

$\blacksquare$

\begin{cor}
\label{Corollary 3.12}Let $X\subseteq L\backslash \{1\}$ and assume that $%
\mathcal{L}$ is an $X$-top lattice. The following conditions are equivalent
for $A\subseteq X:$

\begin{enumerate}
\item $A$ is irreducible;

\item $I(A)$ is strongly irreducible in $(\mathcal{C}(L),\wedge )$;

\item $I(A)$ is irreducible in $(\mathcal{C}(L),\wedge )$.
\end{enumerate}
\end{cor}

\begin{punto}
A \emph{maximal element} in $\mathcal{L}$ is a
maximal element in the poset $(L\backslash \{1\},\leq )$. An element $x\in L$
is called \emph{minimal} in $\mathcal{L}$ iff $x$ is maximal in $\mathcal{L}%
^{0}$. We denote by $Max(\mathcal{L})$ (resp. $Min(\mathcal{L})$) the set of
all maximal (resp. minimal) elements in $\mathcal{L}$. The lattice $\mathcal{%
L}$ is called \emph{coatomic} iff for every element $x\in L\backslash \{1\}$%
, there exists $y\in Max(\mathcal{L})$ such that $x\leq y$. Dually, $%
\mathcal{L}$ is called \emph{atomic} iff for every element $x\in L\backslash
\{0\}$, there exists $y\in Min(\mathcal{L})$ such that $y\leq x$.

Let $A\subseteq L.$ The lattice $\mathcal{L}$ is said to have \textit{the }%
\emph{complete }$A$\emph{-property} iff $\bigwedge\limits_{p\in A\backslash
\{q\}}p\nleq q$ for any $q\in A$. The lattice $\mathcal{L}$ is said to have
the \emph{complete max property} iff $L$ has the complete $Max(\mathcal{L})$%
-property.
\end{punto}

\begin{lem}
\label{Lemma 3.14}Let $\mathcal{L}$ be an $X$-top lattice. If $\mathcal{L}$
is coatomic and $Max(\mathcal{L})\subseteq X,$ then $Max(\mathcal{L}%
)=Max(X). $
\end{lem}

\begin{Beweis}
Let $p\in Max(X)$. Since $\mathcal{L}$ is coatomic, there is $y\in Max(%
\mathcal{L})$ such that $p\leq y$ and so $p=y$ as $Max(\mathcal{L})\subseteq
X$.
\end{Beweis}

$\blacksquare$

The following result recovers and generalizes \cite[Proposition 3.45]%
{Abu2011-a}, \cite[Propositions 5.33, 4.27]{Abu}, and \cite[Proposition 3.40]%
{Abu2011-b}. The additional conditions assumed in these results imply that $%
Max(\mathcal{L})=Max(X)$ (or $Min(\mathcal{L})=Min(X)$ in the dual cases).

\begin{prop}
\label{Proposition 3.15}Let $X\subseteq L\backslash \{1\}.$ The following
are equivalent for $(X,\tau ^{cl}):$

\begin{enumerate}
\item $X$ is $T_{1};$

\item $Max(X)=X=Min(X).$
\end{enumerate}
\end{prop}

\begin{Beweis}
$X$ is $T_{1}$ $\Leftrightarrow $ every singleton is closed $\Leftrightarrow
$ $\{p\}=%
\overline{\{p\}}=V(p)$ $\forall p\in X$ $\Leftrightarrow $ $Max(X)=X=Min(X).$
\end{Beweis}

$\blacksquare$

\begin{thm}
\label{Theorem 3.16 000}Let $X\subseteq L\backslash \{1\}$ and consider $%
(X,\tau ^{cl}).$ Then $Max(X)=Max(\mathcal{C}(\mathcal{L}))$. Moreover, the
following conditions are equivalent:

\begin{enumerate}
\item $X$ is $T_{1}$ and $\mathcal{C}(\mathcal{L})$ satisfies the complete
max property;

\item $X$ is discrete.
\end{enumerate}
\end{thm}

\begin{Beweis}
Let $p\in Max(X)$. Then $p\in \mathcal{C}(L)$ and so $p\in Max(\mathcal{C}(%
\mathcal{L}));$ otherwise, there is $x\in \mathcal{C}(L)\backslash \{1\}$
such that $p\lneq x.$ Since $x\neq 1,$ there is $q\in X$ such that $x\leq q$
and so $p\lneq q$ (a contradiction). For the reverse inclusion, let $x\in
Max(\mathcal{C}(L)).$ Notice that $x=\bigwedge\limits_{p\in A}p$ for some $%
\emptyset \neq A\subseteq X$. Since $A\subseteq \mathcal{C}(L),$ it follows
by the maximality of $x$ in $\mathcal{C}(L)$ that $x=\bigwedge_{p\in A}p=q$
for some $q\in A,$ i.e. $A$ is singleton and $x\in X$. Moreover, $x\in
Max(X) $ as $X\subseteq \mathcal{C}(L)$.

$(1)\Rightarrow (2):$ Assume that $\mathcal{C}(\mathcal{L})$ satisfies the
complete max property. Since $Max(X)=Max(\mathcal{C}(L)),$ we have $%
\bigwedge\limits_{p\in Max(X)\backslash \{q\}}p\nleq q$ for any $q\in Max(X)$%
. Notice that for any $q\in X,$ we have $X=V(\bigwedge\limits_{p\in
X\backslash \{q\}}p)\cup \{q\}$ and by our assumption $q\notin
V(\bigwedge\limits_{p\in X\backslash \{q\}}p)$. Hence, every singleton in $X$
is open, that is $(X,\tau ^{cl})$ is discrete. \bigskip

$(2)\Rightarrow (1):$ Assume that $X$ is discrete and show that $\mathcal{C}(%
\mathcal{L})$ satisfies the complete max property. To show this, suppose
that $q\in X$ and let $Y=X\backslash \{q\}$. Observe that
\begin{equation*}
Y = \overline{Y} = V(I(Y))
\end{equation*}
as $\{q\}$ is open. Hence, $I(Y)\nleq q$, which completes the proof as $%
X=Max(X) = Max(\mathcal{C}(\mathcal{L}))$.
\end{Beweis}

$\blacksquare$

\bigskip

The following result generalizes \cite[Theorem 5.34]{Abu}, \cite[Theorem 4.28%
]{Abu} and \cite[Theorem 3.46]{Abu2011-a}.

\begin{cor}
\label{Theorem 3.16}Let $X\subseteq L\backslash \{1\}$ and consider $(X,\tau
^{cl}).$ If $\mathcal{C}(\mathcal{L})$ satisfies the complete max property,
then the following conditions are equivalent:

\begin{enumerate}
\item $Max(X)=X=Min(X);$

\item $X$ is $T_{2}$;

\item $X$ is $T_{1}$;

\item $X$ is discrete.
\end{enumerate}
\end{cor}

\begin{cor}
\label{Corollary 3.17}Let $X\subseteq L\backslash \{1\}$ and consider $%
(X,\tau ^{cl}).$ Assume that $\mathcal{L}$ satisfies the complete max
property, $\mathcal{L}$ is coatomic and $Max(\mathcal{L})\subseteq \mathcal{C%
}(L).$ The following are equivalent:

\begin{enumerate}
\item $Max(X)=X=Min(X)$;

\item $X$ is $T_{2}$;

\item $X$ is $T_{1}$;

\item $X$ is discrete.
\end{enumerate}
\end{cor}

\begin{Beweis}
Let $p\in Max(\mathcal{C}(\mathcal{L}))$. Since $\mathcal{L}$ is coatomic,
there exists $q\in Max(\mathcal{L})$ such that $p\leq q.$ By assumption, $%
Max(\mathcal{L})\subseteq \mathcal{C}(L)$ whence $p=q.$ So, $Max(\mathcal{L}%
)=Max(\mathcal{C}(\mathcal{L}))$ and the result follows by Corollary {Theorem 3.16}.
\end{Beweis}

$\blacksquare$

A topological space is \emph{regular} \cite{M2000} iff
any non-empty closed set $F$ and any point $x$ that does not belong to $F$
can be separated by disjoint open neighborhoods. A $T_{3}$ space is one
which is both $T_{1}$ and regular. In general, regular spaces need not be
Hausdorff. However, we have a special situation.

\begin{rem}
\label{Proposition 3.19}If $(X,\tau ^{cl})$ is regular, then $(X,\tau ^{cl})$
is $T_{3}.$ To see this, assume that $X$ is regular and let $p\neq q$ be
elements in $X.$ Assume, without loss of generality, that $p\nleq q$ so that
is, $q\notin V(p)$. Since $X$ is regular, there are two disjoint open sets $%
O_{1}$ and $O_{2}$ in $X$ such that $q\in O_{1}$ and $V(p)\subseteq O_{2}$.
Therefore, $X$ is $T_{2}.$
\end{rem}

\index{normal space} A topological space $X$ is \emph{normal} \cite{M2000}
iff any two disjoint closed sets of $X$ can be separated by disjoint open
neighborhoods. The following example shows that the normality of $(X,\tau
^{cl})$ does not guarantee that it is regular.

\begin{ex}
\label{Example 3.20}Let $R$ be a local ring with $|Spec(R)|\geq 2$. Then $%
Spec(R)$ is normal because it has no disjoint non-empty closed sets.
However, $Spec(R)$ is not $T_{1}$ whence not regular by Remark \ref%
{Proposition 3.19}. To see this, notice that the assumption $|Spec(R)|\geq 2$
implies that there is a prime ideal $p$ of $R$ and a maximal ideal $q$ of $R
$ such that $p\lneq q.$ Hence, every open set containing $q$ contains $p$ as
well.
\end{ex}

\subsection{Examples}

Through the rest of this section, $R$ is an associative ring, $M$ is
a non-zero left $R$-module and $LAT(M):=(\mathcal{L}(M),\cap ,+,M,0)$ the
complete lattice of $R$-submodules of $M.$ Moreover, we denote by $Max(M)$
(resp. $\mathcal{S}(M) $) the possibly empty set of maximal (resp. simple) $%
R $-submodules of $M. $ By $L\leq M,$ we mean that $L$ is an $R$-submodule
of $M.$ With abuse of notation, we mean by $I\leq R$ that $I$ is a (two
sided) ideal of $R$.

\begin{punto}
Let $M$ be a left $R$-module. We call an $R$%
-submodule $K\leq M$ :

\emph{prime} \cite{D} iff $K\neq M$ and for any $N\leq M$ and $I\leq R,$ we
have
\begin{equation*}
IN\subseteq K\Rightarrow N\subseteq K%
\text{ or }IM\subseteq K.
\end{equation*}

\emph{first} \cite{AbuC} iff $K\neq 0$ and for any $N\leq K$ and $I\leq R,$
we have
\begin{equation*}
IN=0\Rightarrow N=0\text{ or }IK=0.
\end{equation*}%
\emph{coprime} \cite{Abu} iff $K\neq M$ and for any $I\leq R,$ we have
\begin{equation*}
IM+K=M\text{ or }IM\subseteq K.
\end{equation*}%
\emph{second} \cite{Abu} iff $K\neq 0$ and for any $I\leq R$ we have
\begin{equation*}
IK=K\text{ or }IK=0.
\end{equation*}%
By $Spec^{p}(M)$ (resp. $Spec^{f}(M),$ $Spec^{c}(M),$ $Spec^{s}(M)$) we
denote the spectrum of prime (resp. first, coprime, second) $R$-submodules
of $M.$
\end{punto}

\begin{punto}
\index{fully invariant submodule}
\index{duo module} An $R$-submodule $K$ of $M$ is said to be \emph{fully
invariant} \cite{Abu2011-a} (and we write $L\leq ^{f.i.}M$) iff $%
f(L)\subseteq L$ for all $f\in S:=End(M)$. If every submodule of $M$ is
fully invariant, then $M$ is said to be a \emph{duo module} \cite{Abu2011-a}.
For any $L_{1},L_{2}\leq M,$ we define \begin{equation*}
L_{1}\ast L_{2}=\sum_{f\in Hom(M,L_{2})}f(L_{1})%
\text{\hspace{8pt} and \hspace{8pt}}L_{1}\odot _{M}L_{2}=\bigcap_{f\in
S,f(L_{1})=0}f^{-1}(L_{2});
\end{equation*}%
see \cite{Abu2011-a} and \cite{Abu2011-b}. Notice that if $L_{1}\leq
^{f.i.}M,$ then $L_{1}\ast L_{2}\subseteq L_{1}\cap L_{2}.$
\end{punto}

\begin{punto}
A fully invariant submodule $K\leq ^{f.i.}M$ is:
\emph{fully prime} in $M$ \cite{Abu2011-a} iff $K\neq M$ and for any $%
L_{1},L_{2}\leq ^{f.i.}M,$ we have
\begin{equation*}
L_{1}\ast L_{2}\subseteq K\Rightarrow L_{1}\subseteq K%
\text{ or }L_{2}\subseteq K.
\end{equation*}%
\emph{fully coprime} in $M$ \cite{Abu2011-b} iff $K\neq 0$ and for any $%
L_{1},L_{2}\leq _{R}^{f.i.}M$ we have
\begin{equation*}
K\subseteq L_{1}\odot _{M}L_{2}\Rightarrow K\subseteq L_{1}\text{ or }%
K\subseteq L_{2}.
\end{equation*}%
By $Spec^{fp}(M)$ (resp. $Spec^{fc}(M)$) we denote the spectrum of fully
prime (resp. fully coprime) $R$-submodules of $M.$
\end{punto}

The following example summarizes some facts about some Zariski-like
topologies on several spectra of submodules of a given module.

\begin{ex}
\label{Example 3.2}Consider $X_{1}:=Spec^{p}(M),$ $X_{2}:=Spec^{c}(M),$ $%
X_{3}:=Spec^{fp}(M),$ $X_{4}:=Spec^{s}(M),$ $X_{5}:=Spec^{f}(M)$ and $%
X_{6}:=Spec^{fc}(M).$ Notice that $X_{1},X_{2},X_{3}\subseteq \mathcal{L}%
(M)\backslash \{M\}$ and so one can construct the classical Zariski topology
$\tau _{-}^{cl}$ on any of them as we did for general complete lattices $%
\mathcal{L}=(L,\wedge ,\vee ,1,0)$ and $X\subseteq L\backslash \{1\}.$ On
the other hand, one can construct dual classical Zariski topologies on $\tau
_{-}^{dcl}$ only any of $X_{4},X_{5},X_{6}\subseteq \mathcal{L}(M)\backslash
\{0\}$. Moreover, $M$ is $top^{p}$\emph{-module} (resp. a $top^{c}$\emph{%
-module}, a $top^{fp}$\emph{-module}) if and only if $LAT(M)$ is $X_{1}$-top
(resp. $X_{2}$-top, $X_{3}$-top). On the other hand, $M$ is a $top^{s}$\emph{%
-module} (resp. a $top^{f}$\emph{-module}, a $top^{fc}$-module) iff $LAT(M)$
is dual $X_{4}$-top (resp. dual $X_{5}$-top, dual $X_{6}$-top). The
following table summarize some facts about these spaces.

\newpage

\begin{table}[tbp]
\centering
\begin{adjustbox}{width=1\textwidth}
 \begin{tabular}{ | m{10em} | m{10em}| m{10em} | }
 \hline
   \vspace{5pt} Type \vspace{5pt} & $M\notin Spec^{-}(M)$  & $0\notin Spec^{-}(M)$ \\ \hline
    \vspace{5pt} Subsets of $L$ \vspace{5pt} & \vspace{5pt} $X_1=Spec^{p}(M)$, $X_2=Spec^{c}(M)$, $X_3=Spec^{fp}(M)$ \vspace{5pt} & \vspace{5pt} $X_4=Spec^{f}(M)$, $X_5=Spec^{s}(M)$, $X_6=Spec^{fc}(M)$ \vspace{5pt} \\ \hline
    \vspace{5pt} Variety $V^{-}(N)$ \vspace{5pt} & \vspace{5pt} $\{P\in Spec^{-}(M) \mid N\leq P\}$ \vspace{5pt} & $ \{P\in Spec^{-}(M) \mid P\leq N\}$ \\ \hline
    \vspace{5pt} Subbase $\tau_{-}$ \vspace{5pt} & $\{X\backslash V^{-}(N) \mid N\leq M\}$ & $ \{X\backslash V^{-}(N) \mid N\leq
    M\}$ \\ \hline
    \vspace{5pt} $\tau^{cl}$ or $\tau^{dcl}$ \vspace{5pt} & \vspace{5pt} $\tau_{-}^{cl}$: the topology generated by $\tau_{-}$ \vspace{5pt} &  $\tau_{-}^{dc}$: the topology generated by
    $\tau_{-}$ \\ \hline
    \vspace{5pt} (Dual) $X$-top \vspace{5pt} & \vspace{5pt} $\mathcal{L}$ is  $X_i$-top $\Leftrightarrow \tau_{-}=\tau_{-}^{cl}$ \vspace{5pt} & \vspace{5pt} $\mathcal{L}^0$ is
    $X_j$-top  $\Leftrightarrow \tau_{-}=\tau_{-}^{dcl}$\vspace{5pt} \\ \hline
    \end{tabular}
    \end{adjustbox}
\caption{Examples of spectra of submodules of a given module}
\label{spectra of $M$}
\end{table}
\end{ex}

\begin{ex}
\index{$LAT(M)$} Let $M$ be a local left module over an arbitrary ring $R$,
i.e. $M$ has a unique maximum proper submodule (e.g. the $\mathds{Z}$-module
$\mathds{Z}_{p^{k}}$, $p$ is any prime and $k$ is any positive integer).
Consider $X_{1} = Spec^p(M)$ and $X_{2} = Spec^c(M)$. Then $\mathcal{C}%
^{X_{1}}(LAT(M)) $ and $\mathcal{C}^{X_{2}}(LAT(M))$ satisfy the complete
max property (notice that any maximal submodule is prime and coprime).
\end{ex}

\begin{ex}
\label{Example 3.18}Consider $\mathcal{L}:=LAT(M),$ $X_{1}=Spec^{p}(M)$ and $%
X_{2}=Spec^{c}(M).$ Every maximal submodule of $M$ is a prime and a coprime
submodule, i.e. $Max(M)$ $\subseteq X_{1}$ and $Max(M)\subseteq X_{2}.$ So,
it is enough to assume that $_{R}M$ is a coatomic module satisfying the
complete max property to satisfy the equivalent conditions of Corollary \ref%
{Corollary 3.17}. Moreover, Corollary \ref{Theorem 3.16} applies if $%
RAD^{p}(M):=\mathcal{C}^{X_{1}}(\mathcal{L})$ (resp. $RAD^{c}(M):=\mathcal{C}%
^{X_{2}}(\mathcal{L})$ ) satisfies the complete max property as a lattice.
\end{ex}

\begin{ex}
Consider $\mathcal{L}:=LAT(M),$ $X_{4}=Spec^{f}(M)$ and $X_{5}=Spec^{s}(M).$
Every simple submodule of $M$ is a second and a first submodule of $M,$ i.e.
$\mathcal{S}(M)\subseteq X_{4}$ and $\mathcal{S}(M)\subseteq X_{5}.$ So, it
is enough to assume that $_{R}M$ is an atomic module with the complete min
property to satisfy the equivalent conditions of Corollary \ref{Corollary
3.17} applied to $\mathcal{L}^{0}$. Moreover, Corollary \ref{Theorem 3.16}
applies if $\mathcal{C}^{X_{4}}(\mathcal{L}^{0})$ (resp. $\mathcal{C}%
^{X_{5}}(\mathcal{L}^{0})$) satisfies the complete min property as a lattice.
\end{ex}

\begin{rem}
It was proved in \cite{Abu}, that if $_{R}M$ is a coatomic $top^{c}$-module
satisfying the complete max property, then%
\begin{equation*}
Spec^{c}(M)=Max(M)\Longleftrightarrow X%
\text{ is }T_{2}\Leftrightarrow X\text{ is }T_{1}\Leftrightarrow X\text{ is
discrete.}
\end{equation*}%
A similar result was proved for $Spec^{fp}(M)$ assuming that $_{R}M$ is a
self projective coatomic duo module ($S-PCD$). Notice that it was proved in
\cite[Remark 3.12]{Abu2011-a} that if $_{R}M$ is self projective and duo,
then every maximal submodule is fully prime. Other conditions were assumed
on $M$ in the dual cases to ensure that $\mathcal{S}(M)=Min(X).$ So, Corollary
\ref{Theorem 3.16} generalizes all the corresponding results in \cite%
{Abu2011-a} and \cite{Abu}.
\end{rem}

\subsection{Spectral Spaces}

As before, $\mathcal{L}=(L,\wedge ,\vee ,0,1)$ is a complete lattice.

\begin{punto}
\index{spectral space}
\index{Hochster characterization} \label{Hochster}A topological space is
said to be \emph{spectral} \cite{H1969} iff it is homeomorphic to $Spec(R)$,
the prime spectrum of a commutative ring $R$ with the Zariski topology.
Hochster \cite[Proposition 4]{H1969} characterized such spaces. A
topological space ($X,\tau )$ is spectral if and only if all of the
following four conditions are satisfied:

\begin{enumerate}
\item $X$ is compact;

\item $X$ has a basis of compact open sets closed under finite intersections;

\item $X$ is sober.
\end{enumerate}
\end{punto}

\begin{rem}
\label{Remark 3.21}Let $X\subseteq L\backslash \{1\}.$ If $X$ is finite,
then $(X,\tau ^{cl})$ is spectral:\ By Proposition \ref{Proposition 3.6}, $X$
is $T_{0}$ and sober. The remaining Hochster's conditions in \ref{Hochster}
follow directly from the finiteness of $X$.
\end{rem}

\begin{defn}
\index{radical condition} Let $X\subseteq L\backslash \{1\}$ and consider $%
(X,\tau ^{cl}).$ Set%
\begin{equation}
\mathcal{R}(\mathcal{L}):=\{%
\sqrt{x}\mid x\in L\text{ and }V(x)\text{ is irreducible}\}.
\end{equation}%
We say that $X$ satisfies the \emph{radical condition} iff $\mathcal{R}(%
\mathcal{L})\subseteq X.$
\end{defn}

\begin{lem}
\label{Lemma 3.22}Let $X\subseteq L\backslash \{1\}$ and consider $(X,\tau
^{cl}).$ If $X$ is sober, then $X$ satisfies the radical condition. The
converse holds if $\mathcal{L}$ is $X$-top.
\end{lem}

\begin{Beweis}
Let $X$ be sober. Let $x\in \mathcal{R}(\mathcal{L})$. Since $X$ is sober, $%
V(x)=\overline{\{p\}}\overset{\text{Proposition \ref{Proposition 3.4}}}{=}V(p)$ for some
$p\in X.$ It follows by Lemma \ref{lemma 1.19} (2) that $\sqrt{x}=p\in X.$

For the converse, assume that $\mathcal{L}$ is $X$-top. Let $F$ be a closed
irreducible subset of $X.$ Since $\mathcal{L}$ is $X$-top, $F=V(x)$ for some
$x\in L$. By our hypothesis, $\sqrt{x}\in X.$ By Lemma \ref{lemma 1.19} (3),
$F=V(x)=V(\sqrt{x})$ and so $\sqrt{x}$ is the unique generic point of $F$
(the uniqueness is obvious). Therefore, $X$ is sober.
\end{Beweis}

$\blacksquare$

\begin{prop}
\label{Proposition 3.23}Let $X\subseteq L\backslash \{1\}$ and assume that $%
\mathcal{L}$ is an $X$-top lattice. If $\mathcal{L}$ satisfies the ACC, then
every subset of $(X,\tau )$ is compact.
\end{prop}

\begin{Beweis}
Let $A\subseteq X$ and suppose that $\mathcal{O}=\{X\backslash V(x_{i}))\mid
x_{i}\in L,$ $i\in I\}$ is an open cover for $A$. Since $\mathcal{L}$
satisfies the ACC, $\bigvee_{i\in I}x_{i}=\bigvee_{j\in J}x_{j}$ for some
finite subset $J\subseteq I$. Notice that
\begin{equation*}
A\subseteq \bigcup\limits_{i\in I}(X\backslash V(x_{i}))=X\backslash
V(\bigvee_{i\in I}x_{i}))=X\backslash V(\bigvee_{i\in
I}x_{i}))=\bigcup\limits_{i\in J}(X\backslash V(x_{j})),
\end{equation*}%
i.e. $\{X\backslash V(x_{j})\mid j\in J\}$ is a finite subcover of $\mathcal{O}
$ for $A$.
\end{Beweis}

$\blacksquare$

\begin{prop}
\index{$\mathcal{C}(\mathcal{L})$} \label{Proposition 3.24}Let $X\subseteq
L\backslash \{1\}$ and assume that $\mathcal{L}$ is an $X$-top lattice. The
following conditions are equivalent:

\begin{enumerate}
\item $\mathcal{C}(\mathcal{L})$ satisfies the ACC;

\item Every subset of $(X,\tau )$ is compact;

\item Every open set in $(X,\tau )$ is compact.
\end{enumerate}
\end{prop}

\begin{Beweis}
$(1\Rightarrow 2):$ Consider the complete lattice $(\mathcal{C}(L),\wedge ,%
\tilde{\vee},\sqrt{0},1).$ Since $V(x)=V(\sqrt{x})$ for every $x\in L,$ we
conclude that $\mathcal{C}:=(\mathcal{C}(L),\wedge ,\tilde{\vee},\sqrt{0},1)$
is an $X$-top lattice. By our assumption, $\mathcal{C}$ satisfies the ACC
and so every subset of $X$ is compact by Proposition \ref{Proposition 3.23}.

$(3\Rightarrow 1):$ Let $a_{1}\leq a_{2}\leq a_{3}\leq \cdots $ be an
ascending chain in $\mathcal{C}(\mathcal{L})$. Notice that $X\backslash
V(a_{1})\subseteq X\backslash V(a_{2})\subseteq X\backslash
V(a_{3})\subseteq \cdots $. Setting $b=\tilde{\bigvee }_{i=1}^{\infty }a_{i}$%
, we observe that
\begin{equation*}
X\backslash V(b)=X\backslash V(\bigvee_{i=1}^{\infty }a_{i})=X\backslash
\bigcap_{i=1}^{\infty }V(a_{i})=\bigcup_{i=1}^{\infty }(X\backslash
V(a_{i})).
\end{equation*}%
By our assumption, the open set $X\backslash V(b)$ is compact and so $%
X\backslash V(b)=\bigcup\limits_{i=1}^{n}X\backslash V(a_{i})=X\backslash
V(a_{n})$ for some $n\in \mathbb{N},$ i.e. $b=a_{n}$ and the ascending chain
stabilizes.
\end{Beweis}

$\blacksquare$

\begin{cor}
\label{Corollary 3.25}Let $X\subseteq L\backslash \{1\}$ and $\mathcal{L}$
be an $X$-top lattice such that $\mathcal{C}(\mathcal{L})$ satisfies the
ACC. Then $(X,\tau )$ is spectral $\Leftrightarrow $ $(X,\tau )$ is sober.
\end{cor}

\begin{Beweis}
By Proposition \ref{Proposition 3.6} $X$ is $T_{0}.$ The result follows now
using Proposition \ref{Proposition 3.24} and Hochster's characterization for
spectral spaces \ref{Hochster}.
\end{Beweis}

$\blacksquare$

\noindent In \cite{BH2008-b}, the so called \emph{finer patch topology} was
used to prove that for any left module $M$ over an associative ring $R,$ and
$X=Spec^{p}(M),$ the classical Zariski topology $(X,\tau ^{cl})$ is a
spectral space provided that the ACC holds for intersections of prime submodules
of $M$.

\begin{punto}
\label{finer}Let $X\subseteq L\backslash \{1\}.$ The \textit{finer patch
topology%
\index{finer patch topology} }$\tau ^{\mathfrak{fp}}$ on $X$ is the one
whose subbase is%
\begin{equation}
\mathcal{B}=\{V(x)\cap X\backslash V(y)\mid x,y\in L\}.
\end{equation}%
It is clear that $\tau ^{cl}\subseteq \tau ^{\mathfrak{fp}}$. So, if $\tau ^{%
\mathfrak{fp}}$ is compact, then $\tau ^{cl}$ is compact.
\end{punto}

\begin{ex}
Let $\mathcal{P}$ be the set of all prime numbers in $\mathds{Z}$ and
consider the ring $R = \prod_{p\in \mathcal{P}}\mathds{Z}_p.$ Then the finer
patch topology associated with $Spec(R)$ is not compact while, trivially,
the classical Zariski topology is compact. In general, if $R$ is a ring with
zero dimension and $Spec(R)$ is infinite, then the finer patch topology
associated with $Spec(R)$ is not compact while, trivially, the classical
Zariski topology is compact.
\end{ex}

\begin{prop}
\index{finer patch topology} \label{Proposition 3.26}Let $X\subseteq
L\backslash \{1\}.$ If $(X,\tau ^{\mathfrak{fp}})$ is compact, then $(X,\tau
^{cl})$ is spectral.
\end{prop}

\begin{Beweis}
Assume that $(X,\tau ^{\mathfrak{fp}})$ is compact. We apply Hochster's
characterizations of spectral spaces to prove that $(X,\tau ^{cl})$ is
spectral. Notice that $(X,\tau ^{cl})$ is $T_{0}$ by Proposition \ref%
{Proposition 3.6} and is compact since $\tau ^{cl}\subseteq \tau ^{\mathfrak{%
fp}}.$

\textbf{Claim I:}$\ (X,\tau ^{cl})$ is sober.

Let $Y\subseteq X$ be a closed irreducible set in $(X,\tau ^{cl})$. Then $Y%
\overset{%
\text{Proposition \ref{Proposition 3.4}}}{=}\bigcup\limits_{p\in Y}V(p).$ On the other
hand $Y$ is closed in $(X,\tau ^{\mathfrak{fp}}),$ whence compact in $%
(X,\tau ^{\mathfrak{fp}})$ (recall that every closed subset of a compact
space is compact). Therefore, the open cover $\mathcal{O}:=\{V(p):p\in Y\}$
of $Y$ has a finite subcover $\{V(p_{1}),\cdots ,V(p_{n})\},$ i.e. $%
Y=\bigcup\limits_{i=1}^{i=n}p_{i}$. But $Y$ is irreducible, whence $%
Y=V(p_{k})$ for some $k\in \{1,2,\cdots ,n\}.$ Clearly, $p_{k}$ is the
unique generic point of $Y.$

\textbf{Claim II: }$X$ has a basis of compact open sets closed under finite
intersections.

We prove this claim in two steps.

\textbf{Step 1:}\ Every basic open subset of $(X,\tau ^{cl})$ is compact.

Let $B$ be a basic open subset of $(X,\tau ^{cl})$, i.e. $%
B=\bigcap\limits_{i=1}^{i=n}X\backslash V(x_{i})$ for some $\{x_{1},\cdots
,x_{n}\}\subseteq L.$ Observe that $X\backslash V(x_{i})$ is closed in $%
(X,\tau ^{\mathfrak{fp}})$ $\forall i\in \{1,2,\cdots ,n\}$. So, $B$ is
closed in ($X,\tau )$, whence compact in $(X,\tau ^{\mathfrak{fp}}).$ Since $%
\tau ^{cl}\subseteq \tau ^{\mathfrak{fp}},$ $B$ is compact in $(X,\tau
^{cl}).$

\textbf{Step 2:} The collection of open compact subsets of $(X,\tau ^{cl})$
is closed under arbitrary intersections.

Let $U$ be an open compact set in $(X,\tau ^{cl}).$ Then we can write $%
U=\bigcup\limits_{i=1}^{n}\bigcap\limits_{j=1}^{m_{i}}X\backslash V(x_{ij})$
for some subset $\{x_{ij}\mid j=1,2,\cdots ,m_{i},$ $i=1,\cdots ,n\}$ (the
union is finite because of the compactness of $U$). Notice that $U$ is
closed in $(X,\tau ^{\mathfrak{fp}}).$ So, any intersection of open compact
subsets in $(X,\tau ^{cl})$ is closed in $(X,\tau ^{\mathfrak{fp}});$ so it
is compact in $(X,\tau ^{\mathfrak{fp}}),$ whence compact in $(X,\tau ^{cl})$%
.
\end{Beweis}

$\blacksquare$

\begin{ex}
\label{Example 3.27}The ring of integers $\mathbb{Z}$ is Noetherian and so
the finer patch topology on $Spec(\mathbb{Z})$ is compact because the ACC is
satisfied on the radical ideals by \cite[Theorem 2.2]{BH2008-b}. This
example shows that $(X,\tau ^{\mathfrak{fp}})$ can be compact although $X$
is infinite.
\end{ex}

\begin{ex}
\label{Example 3.28}Let $L$ be infinite and be such that the elements of $%
X:=L\backslash \{0,1\}$ are not comparable (notice that for all $a\neq b$ in
$X$ we have $a\wedge b=0$ and $a\vee b=1$). Notice that $(X,\tau ^{\mathfrak{%
fp}})$ is not compact, whereas $(X,\tau ^{cl})$ is compact because it is the
cofinite topology on $X.$ Notice also that $\mathcal{C}(\mathcal{L})$
satisfies the ACC and every element in $\mathcal{C}(\mathcal{L})$ can be
written as an irredundant meet of elements in $X$, but this guarantees the
compactness for the finer patch topology. Observe that $\mathcal{L}$ is not $%
X$-top and $(X,\tau ^{cl})$ is not sober and hence not spectral.
\end{ex}

\begin{prop}
\index{irreducible subset} \label{Proposition 4.11}Let $X\subseteq
L\backslash \{1\}$ and consider $(X,\tau ^{cl})$. If $V(x)$ is reducible for
some $x\in L$, then $V(x)=\bigcup\limits_{i=1}^{n}V(x_{i})$ for some
elements $x_{1},x_{2},\cdots ,x_{n}\in L$, where $V(x_{i})\subsetneqq V(x)$
for all $i=1,2,\cdots ,n$.
\end{prop}

\begin{Beweis}
Let $V(x)$ be reducible for some $x\in L$, i.e. $V(x)=F_{1}\cup F_{2}$ where
both $F_{1}$ and $F_{2}$ are closed proper subsets of $V(x)$. Suppose that $%
F_{1}=\bigcap_{i\in I}\bigcup_{j=1}^{n_{i}}V(x_{ij})$ and $%
F_{2}=\bigcap_{l\in L}\bigcup_{k=1}^{m_{l}}V(y_{lk})$ for some $%
\{x_{ij}\},\{y_{lk}\}\subseteq L.$ Since $F_{1}$ and $F_{2}$ are proper
subsets of $V(x),$ we have $V(x)\nsubseteq \bigcup_{j=1}^{n_{i_{0}}}V(x_{{%
i_{0}}j})$ for some $i_{0}\in I$ and $V(x)\nsubseteq
\bigcup_{k=1}^{m_{l_{0}}}V(y_{{l_{0}}k})$ for some $l_{0}\in L$, whence $%
V(x)\neq \bigcup_{j=1}^{n_{i_{0}}}V(x_{{i_{0}}j})\cap V(x)$ and $V(x)\neq
\bigcup_{k=1}^{m_{l_{0}}}V(y_{{l_{0}}k})\cap V(x)$. Set $x_{r}:=x_{{i_{0}}%
r}\vee x$ for $r=1,2,\cdots ,n_{i_{0}}$ and $x_{n_{i_{0}}+r}=y_{{l_{0}}%
r}\vee x$ for $r=1,2,\cdots ,m_{l_{0}}$ and let $n:=n_{i_{0}}+m_{l_{0}}$. By
construction, $V(x)=\bigcup_{r=1}^{n}V(x_{r})$, where each $V(x_{r})$ is a
proper subset of $V(x)$.
\end{Beweis}

$\blacksquare$

As a direct consequence of Proposition \ref{Proposition 4.11}, we obtain the
following result which recovers \cite[Proposition 2.26]{BH2008-a} proved for
the prime spectrum of a module over a ring.

\begin{cor}
\label{Corollary 3.32}Let $X\subseteq L\backslash \{1\}$ and assume that $%
\left\vert X\right\vert \geq 2.$ If $(X,\tau ^{cl})$ is $T_{2},$ then there
exist $x_{1},x_{2},\cdots ,x_{n}\in L$ such that $X=\bigcup%
\limits_{i=1}^{n}V(x_{i})$, while $X\neq V(x_{i})$ for all $i=1,2,\cdots ,n$.
\end{cor}

The radical condition is automatically satisfied by the spectrum of prime
submodules of a given left module over an associative ring by \cite[Theorem
3.4, Corollary 3.6]{BH2008-a}. However, we need to check it when dealing
with other cases.

The following result generalizes \cite[Theorem 3.2]{BH2008-b}:

\begin{thm}
\index{$\mathcal{C}(\mathcal{L})$}
\index{finer patch topology}
\index{radical condition} \label{Theorem 3.33}Let $X\subseteq L\backslash
\{1\}$ and consider $(X,\tau ^{cl})$. If $X$ satisfies the radical
condition, then $\mathcal{C}(\mathcal{L})$ satisfies the ACC if and only if $%
(X,\tau ^{\mathfrak{fp}})$ is compact. It follows that, If $\mathcal{C}(%
\mathcal{L})$ satisfies the ACC and $X$ satisfies the radical condition,
then $(X,\tau ^{cl})$ is spectral.
\end{thm}

\begin{Beweis}
Assume that $\mathcal{C}(\mathcal{L})$ satisfies the ACC and that $X$
satisfies the radical condition. We need only to prove that $(X,\tau ^{%
\mathfrak{fp}})$ is compact since it will follow then, by Proposition \ref%
{Proposition 3.26}, that $(X,\tau ^{cl})$ is spectral.

Suppose that $(X,\tau ^{\mathfrak{fp}})$ is not compact, i.e. there is an
open cover $\mathcal{A}$ in $\tau ^{\mathfrak{fp}}$ for $X$ which does not
have a finite subcover for $X.$

Let
\begin{equation*}
\mathbb{E}:=\{x\in \mathcal{C}(L)\mid V(x)%
\text{ is not covered by a finite subcover of }\mathcal{A}\}.
\end{equation*}%
Observe that $\sqrt{0}\in \mathbb{E},$ i.e. $\mathbb{E}\neq \emptyset .$
Since $\mathcal{C}(\mathcal{L})$ satisfies the ACC, $\mathbb{E}$ has a
maximal element $p$. Notice that $V(p)\neq \emptyset .$

\textbf{Case 1: }$p\notin X.$ Since $X$ satisfies the radical condition, $%
V(p)$ is reducible and it follows by Proposition \ref{Proposition 4.11} that
$V(p)=\bigcup\limits_{i=1}^{n}V(x_{i})$ for some $x_{1},\cdots ,x_{n}\in
\mathcal{C}(L)$ (see Lemma \ref{lemma 1.19} (3)), where $V(x_{i})\subsetneqq
V(p),$ whence $p\lneq x_{i},$ for all $i\in \{1,\cdots ,n\}$. Since $p$ is
maximal in $\mathbb{E}$ and $p\lneq x_{i},$ $V(x_{i})$ is covered by a
finite subcover of $\mathcal{A}$ for all $i\in \{1,\cdots ,n\}$. Hence $V(p)$
is covered by a finite subcover of $\mathcal{A}$, a contradiction.

\textbf{Case 2: }$p\in X.$ It follows that $p\in O$ for some $O\in \mathcal{A%
}$ and so $p\in B$, where $B$ is a basic open subset of $O$. Assume that
\begin{equation*}
B=\bigcap\limits_{i=1}^{n}(V(x_{i})\cap X\backslash V(y_{i}))\text{ for some
subset }\{x_{1},\cdots ,x_{n},y_{1},\cdots ,x_{n}\}\subseteq L.
\end{equation*}
Observe that $z_{i}:=y_{i}\vee p\nleq p$ as $y_{i}\nleq p$ $\forall $ $i\in
\{1,2,\cdots n\}$.

\textbf{Claim:} $V(p)\cap \bigcap\limits_{i=1}^{n}X\backslash
V(z_{i})\subseteq B$. To prove this claim, let $q\in V(p)\cap
\bigcap\limits_{i=1}^{n}X\backslash V(z_{i})$ for all $i\in \{1,2,\cdots
,n\} $, whence $p\leq q$ and $y_{i}\vee p\nleq q$ for all $i\in \{1,2,\cdots
,n\}. $ It follows that $p\leq q$ and $y_{i}\nleq q$ for all $i\in
\{1,2,\cdots ,n\}.$ But $x_{i}\leq p$ for all $i\in \{1,2,\cdots ,n\}$
whence $x_{i}\leq q $ and $y_{i}\nleq q$ for all $i\in \{1,2,\cdots ,n\}$,
i.e. $q\in \bigcap\limits_{i=1}^{n}(V(x_{i})\cap X\backslash V(y_{i}))=B$ as
claimed.

Now, notice that for all $i\in \{1,2,\cdots ,n\}$, we have $p\lneq z_{i}$
and so $V(z_{i})$ is covered by a finite subcover $\mathcal{A}_{i}$ of $%
\mathcal{A}.$ On the other hand, $V(p)\backslash
\bigcup\limits_{i=1}^{n}V(z_{i})=V(p)\cap
\bigcap\limits_{i=1}^{n}X\backslash V(z_{i})\subseteq B\subseteq O$. Hence $%
\{O\}\cup \mathcal{A}_{1}\cup \mathcal{A}_{2}\cup \cdots \cup \mathcal{A}%
_{n} $ is a finite subcover of $\mathcal{A}$ for $V(p)$, which is a
contradiction.

Therefore, $(X,\tau ^{\mathfrak{fp}})$ is compact.

Conversely, assume that $(X,\tau ^{\mathfrak{fp}})$ is compact. Suppose that
$\mathcal{C}(\mathcal{L})$ does not satisfy the ACC. Then there is an
infinite strictly increasing chain $a_{1}\lneq a_{2} \lneq ....$ of elements
in $\mathcal{C}(\mathcal{L})$. Since $(X,\tau ^{\mathfrak{fp}})$ is compact,
$V(a_{1})$ is compact as it is closed. But one can check that $%
\{V(a_{i})\cap(X\backslash V(a_{ i+1}))\mid i=1, 2, ...\}\cup
\{V(\bigvee\limits_{i=1}\limits^{\infty}a_i)\}$ is an open cover for $%
V(x_{1})$ which does not have a finite subcover, a contradiction.
\end{Beweis}

$\blacksquare$

\begin{rem}
\label{Remark 3.43}Let $X\subseteq L\backslash \{1\}.$ The radical condition
in Theorem \ref{Theorem 3.33} is necessary for $(X,\tau ^{cl})$ to be
spectral. Recall that this condition is satisfied if $X$ is sober (see Lemma %
\ref{Lemma 3.22}).
\end{rem}

\begin{defn}
Let $X\subseteq L\backslash \{1\}.$ An element $p\in X$ is called \emph{%
minimal in }$X$\emph{\ over }$x\in L$ iff $p=q$ whenever $x\leq q\leq p$ for
some $q\in X$.
\end{defn}

\begin{cor}
\label{Corollary 3.34}Let $X\subseteq L\backslash \{1\}.$ Assume that $%
\mathcal{C}(\mathcal{L})$ satisfies the ACC, and that for any $x\in \mathcal{%
C}(L)\backslash (X\cup \{1\})$ with $V(x)\neq \emptyset $ there is a
completely strongly irreducible minimal element in $X$ over it with respect
to $(\mathcal{C}(L),\wedge )$. Then $(X,\tau ^{\mathfrak{fp}})$ is compact
(and consequently $(X,\tau ^{cl})$ is spectral).
\end{cor}

\begin{Beweis}
We claim that $X$ satisfies the radical condition. Let $x\in \mathcal{R}(%
\mathcal{L})\backslash X.$ In particular, $V(x)\neq \emptyset .$ Let $p$ be
a completely strongly irreducible minimal element in $X$ over $x$. Then $%
\bigwedge\limits_{q\in V(x)\backslash \{p\}}q\nleq p$ (otherwise, $%
\bigwedge_{q\in V(x)\backslash \{p\}}q\leq p$ and the complete strong
irreducibility of $p$ would imply that $q\lneq p$ for some $q\in V(x)$
contradicting the minimality of $p$ over $x$). Therefore, $%
V(x)=V(\bigwedge_{q\in V(x)\backslash \{p\}}q)\cup V(p)$ a union of proper
closed subsets and so $V(x)$ is reducible, a contradiction. So, $X$
satisfies the radical condition. Now, the hypotheses of Theorem \ref{Theorem
3.33} are satisfied and it follows that $(X,\tau ^{\mathfrak{fp}})$ is
compact and consequently $(X,\tau ^{cl})$ is spectral.
\end{Beweis}

$\blacksquare$

\begin{thm}
\index{finer patch topology} \label{Theorem 3.36}Let $X\subseteq L\backslash
\{1\}$ and consider $(X,\tau ^{cl})$. Assume that $\mathcal{C}(\mathcal{L})$
satisfies the DCC and that $Min(X)\subseteq SI(\mathcal{C}(\mathcal{L}))$.
Then $(X,\tau ^{\mathfrak{fp}})$ is compact if and only if $V(p)$ is finite $%
\forall p\in Min(X)$.
\end{thm}

\begin{Beweis}
Assume that $\mathcal{C}(\mathcal{L})$ satisfies the DCC and that $%
Min(X)\subseteq SI(\mathcal{C}(\mathcal{L}))$. We show that $(X,\tau ^{%
\mathfrak{fp}})$ is compact. Notice first of all that $X=\bigcup\limits_{p%
\in Min(X)}V(p)$, since $\mathcal{C}(\mathcal{L})$ satisfies the DCC.

\textbf{Claim:}\ $Min(X)$ is finite. To prove this claim, notice that $%
\bigwedge\limits_{p\in Min(X)}p=\bigwedge\limits_{i=1}^{n}p_{i}$ for some $%
\{p_{1},p_{2},\cdots ,p_{n}\}\subseteq Min(X)$ (since $\mathcal{C}(\mathcal{L%
})$ satisfies the DCC). So, $\bigwedge_{i=1}^{n}p_{i}\leq p$ for all $p\in
Min(X)$. By assumption, $Min(X)\subseteq SI(\mathcal{C}(\mathcal{L}))$,
whence $p=p_{i}$ for some $i\in \{1,2,\cdots ,n\}.$ Consequently, $Min(X)$
is finite.

If $V(p)$ is finite $\forall p\in Min(X)$, then $X$ is finite, whence $%
(X,\tau ^{\mathfrak{fp}})$ is compact.

Conversely, suppose that $(X,\tau ^{\mathfrak{fp}})$ is compact and that $%
V(p)$ is infinite for some $p\in Min(X)$.

\textbf{Case 1:} $V(p)$ contains an infinite chain $p=x_{1}\lneq x_{2}\lneq
....$ which does not stabilize. Consider the open cover $\mathcal{A}
:=\{V(x_{i})\cap (X\backslash V(x_{i+1})\mid i=1, 2, ..\}\cup
\{V(\bigvee\limits_{i=1}\limits^{\infty}x_i)\}$ for $V(p)$. Clearly $%
\mathcal{A}$ has no finite subcover for $V(x_{1})$, whence $(X,\tau ^{%
\mathfrak{fp}})$ is not compact, a contradiction.

\textbf{Case 2:} $V(p)$ does not contain any infinite chain. It follows that
there is an infinite subset $A\subseteq V(p)$ of incomparable elements.
Since $\mathcal{C}(\mathcal{L})$ satisfied the DCC, it follows that $%
\bigwedge\limits_{x\in A}x = \bigwedge\limits_{x\in F}x$ for some finite
subset $F\subseteq A$. Since $A$ is infinite, there is $q\in A\backslash F$
such that $p\lneq q$ for some $p\in F$, a contradiction.
\end{Beweis}

$\blacksquare$

\begin{lem}
\index{$X-top$} \label{DCC implies SIX} Let $X\subseteq L\backslash \{1\}$
and $\mathcal{L}$ be an $X$-top lattice. Assume that $\mathcal{C}(\mathcal{L}%
)$ satisfies the DCC. Then $X\subseteq SI(\mathcal{C}(\mathcal{L}))$.
\end{lem}

\begin{Beweis}
Since $\mathcal{L}$ is an $X$-top lattice, we have $\tau =\tau ^{cl}$.
Notice that for every $p\in X$, the singleton $\{p\}$ is irreducible in $%
(X,\tau )$, whence $p=I(\{p\})$ is strongly irreducible in $(\mathcal{C}%
(L),\wedge )$ by Proposition \ref{Proposition 3.11}.
\end{Beweis}

$\blacksquare$

\begin{cor}
\label{Corollary 3.37}Let $X\subseteq L\backslash \{1\}$ and $\mathcal{L}$
be an $X$-top lattice. If $\mathcal{C}(\mathcal{L})$ satisfies the DCC, then
$(X,\tau ^{\mathfrak{fp}})$ is compact if and only if $V(p)$ is finite $%
\forall p\in Min(X)$.
\end{cor}

\begin{Beweis}
Follows directly by applying Lemma \ref{DCC implies SIX} and Theorem \ref%
{Theorem 3.36}.
\end{Beweis}

$\blacksquare$

\begin{ex}
\label{Example 3.41}Let $\mathcal{L}=(L,\wedge ,\vee ,1,0)$ be a complete
lattice, where $L$ is an infinite ascending chain $x_{1}\leq x_{2}\leq
\cdots $ endowed with a maximum element $1$ such that $\bigvee\limits_{i\in
I}x_{i}=1$ for every infinite subset $I\subseteq \mathbb{N}$. Let $%
X=L\backslash \{1\}$. Then $\mathcal{C}(\mathcal{L})$ satisfies the DCC, and
$Min(X)\subseteq SI(\mathcal{C}(\mathcal{L})$. Hence, $\tau ^{\mathfrak{fp}}$
is not compact by Theorem \ref{Theorem 3.36} because $V(x_{1})$ is infinite.
Moreover, every descending chain of $(X, \tau)$ is a spectral subspace.
\end{ex}

In what follows, $R$ is a ring, $M$ is a left $R$-module and consider $%
\mathcal{L}:=LAT(M),$ the complete lattice of left $R$-submodules of $M.$

\begin{ex}
\label{Example 3.38}Let $X=Spec^{p}(M),$ the spectrum of prime $R$%
-submodules of $M.$ By \cite[Theorem 3.4 (i)]{BH2008-a}, $Spec^{p}(M)$
satisfies the radical condition. Therefore, Theorem \ref{Theorem 3.33}
recovers \cite[Theorem 3.2]{BH2008-b} as a special case.
\end{ex}

\begin{ex}
\label{Example 3.39}Let $_{R}M$ be Noetherian and $X=SI(M),$ the spectrum of
strongly irreducible $R$-submodules of $M,$ whence $\mathcal{L}$ is $X$-top.
By \cite[Proposition 2.7]{AbuC}, $SI(M)$ satisfies the radical condition.
Therefore $(SI(M),\tau ^{\mathfrak{fp}})$ is compact and $(SI(M),\tau )$ is
spectral.
\end{ex}

\begin{ex}
\label{Example 3.42}Applying Theorem \ref{Theorem 3.33}, we obtain several
examples of spectral spaces:

\begin{enumerate}
\item If $_{R}M$ is duo and $\mathcal{C}(\mathcal{L})$ satisfies the ACC,
then $Spec^{fp}(M)$ is spectral (notice that $Spec^{fp}(M)$ satisfies the
radical condition by \cite[Proposition 3.30]{Abu2011-a}).

\item If $_{R}M$ is duo and $\mathcal{H}(L)$ satisfies the DCC, then $%
X=Spec^{fc}(M)$ is spectral (notice that $X=Spec^{fc}(M)$ satisfies the
radical condition by \cite[Proposition 3.25]{Abu2011-b}).

\item If $_{R}M$ is a completely distributive $top^{c}$-module and $\mathcal{%
C}(\mathcal{L})$ satisfies the ACC, then $Spec^{c}(M)$ is spectral (notice
that $X=Spec^{c}(M)$ satisfies the radical condition by \cite[Proposition
5.19 (i)]{Abu}).

\item If $_{R}M$ is a $top^{s}$-module and $\mathcal{H}(L)$ satisfies the
DCC, then $Spec^{s}(M)$ is spectral (notice that $X=Spec^{s}(M)$ satisfies
the radical condition by \cite[Proposition 4.14 (i)]{Abu}).

\item If $_{R}M$ is a $top^{f}$-module, $I(A)$ is first for every
irreducible subset $A\subseteq Spec^{f}(M)$ and $\mathcal{H}(L)$ satisfies
the DCC, then $Spec^{f}(M)$ is spectral (notice that the assumption on the
irreducible subsets of $X=Spec^{f}(M)$ is equivalent to $X$ satisfying the
radical condition by \cite[Remark 4.25]{AbuC}).
\end{enumerate}
\end{ex}

\section{\label{alg-vs-topolog}Algebraic versus Topological Properties}

As before, $\mathcal{L}=(L,\wedge ,\vee ,0,1)$ is a complete lattice. In
this section we study the interplay between the topological properties of $%
(X,\tau ^{cl})$ where $X\subseteq L\backslash \{1\}$ (or $(X,\tau ^{dcl})$
where $X\subseteq L\backslash \{0\}$) and the algebraic properties of $%
\mathcal{L}$. Applications will be given to the special case $\mathcal{L}%
=LAT(_{R}M),$ where $R$ is a ring and $M$ is a left $R$-module.

\begin{punto}
We say that an element $x\in L$ is%
\index{finitely constructed} is \emph{finitely constructed} in $\mathcal{L}$
iff $x$ cannot be written as an infinite irredundant join of elements of $L$%
, that is, for any collection $\{x_{i}\}_{i\in I}\subseteq L$ such that $%
\bigvee\limits_{i\in I}x_{i}=x$, there is a finite sub-collection $%
\{x_{j}\}_{j\in J}$ of $\{x_{i}\}_{i\in I}$ with $\bigvee\limits_{j\in
J}x_{j}=x$. An element $x$ is called \emph{countably finitely constructed}
\index{countably finitely constructed} in $\mathcal{L}$ iff $x$ cannot be
written as an infinite countable irredundant join of elements of $L$, i.e.
for any countable collection $\{x_{i}\}_{i\in I}\subseteq L$ with $%
\bigvee\limits_{i\in I}x_{i}=x,$ there is a finite sub-collection $%
\{x_{j}\}_{j\in J}$ of $\{x_{i}\}_{i\in I}$ with $\bigvee\limits_{j\in
J}x_{j}=x$. An element $x$ is called \emph{countably constructed}%
\index{countably constructed} in $\mathcal{L}$ iff $x$ cannot be written as
an uncountable irredundant join of elements of $L$.
\end{punto}

We collect first some remarks:

\begin{rems}
\index{$X-top$} \label{Remark 3.44}Let $\mathcal{L}$ be an $X$-top lattice, $%
X\subseteq L\backslash \{1\}$ and consider the topological space $(X,\tau )$.

\begin{enumerate}
\item The following are equivalent:

\begin{enumerate}
\item $(X,\tau )$ is irreducible;

\item $%
\sqrt{0}\in SI(\mathcal{C}(\mathcal{L}));$

\item If $X=\bigcup\limits_{i\in I}V(x_{i})$, then either $I$ is infinite or
there is $i_{0}\in I$ such that $x_{i_{0}}$ is a lower bound for $X$.
\end{enumerate}

\item $(X,\tau )$ is $T_{1}$ if and only if $Max(X)=X$.

\item $(X,\tau )$ is Noetherian $\Leftrightarrow $ $\mathcal{C}(\mathcal{L})$
satisfies the ACC $\Leftrightarrow $ each set in $X$ is compact $%
\Leftrightarrow $ each open set in $X$ is compact.

\item $(X,\tau )$ is Artinian $\Leftrightarrow $ $\mathcal{C}(\mathcal{L})$
satisfies the DCC $\Leftrightarrow $ every closed cover for any subset of $X$
has a finite subcover.

\item $(X,\tau )$ is (countably) compact if and only if $1$ is (countably)
finitely constructed in $\mathcal{C}(\mathcal{L})$.

\item If $SI(\mathcal{C}(\mathcal{L}))\subseteq X$, then $(X,\tau )$ is
sober.

\item If $X$ satisfies the radical condition, then $(X,\tau )$ is sober.

\item Assume that $\mathcal{C}(\mathcal{L})$ satisfies the complete max
property. Then, $(X,\tau )$ is $T_{1}$ $\Leftrightarrow $ $(X,\tau )$ is
discrete.

\item If $(X,\tau )$ is an atomic, Lindelof (compact) and $V(p)$ is open $%
\forall $ $p\in Min(X)$, then $Min(X)$ is countable (finite).

\item $V(x)$ is irreducible for every $x\in X.$
\end{enumerate}
\end{rems}

\begin{Beweis}
Let $\mathcal{L}$ be an $X$-top lattice.

\begin{enumerate}
\item $(a\Leftrightarrow b)$ Apply Corollary \ref{Corollary 3.12} to $%
V(0)=X. $

$(a\Rightarrow c)$ Suppose that $X=\bigcup\limits_{i\in I}V(x_{i})$ with $I$
finite. Since $X$ is irreducible, $V(x_{i_{0}})=X$ for some $i_{0}\in I$
whence $x_{i_{0}}$ a lower bound for $X$.

$(c\Rightarrow a)$ Suppose that $X=V(x)\cup V(y)$ for some $x,y\in L$. By
our assumption, $x$ is a lower bound for $X$ whence $X=V(x)$ or $y$ is a
lower $X$ whence $X=V(y)$. Therefore, $X$ is irreducible.

\item Apply Proposition \ref{Proposition 3.15} to $(X,\tau )=(X,\tau ^{cl})$.

\item It is easy to check that the first two statements are equivalent. The
remaining equivalences follow by applying Proposition \ref{Proposition 3.4}
to $(X,\tau )=(X,\tau ^{cl})$.

\item Notice that any open set in $X$ has the form $X\backslash V(x)$ where $%
x\in \mathcal{C}(L)$. The equivalence of the first two statements is
straightforward. We claim that they are equivalent to the third statement.

Assume that $\mathcal{C}(\mathcal{L})$ satisfies the DCC. Let $U\subseteq X$
and $\{V(x)\mid x\in A\}$ be a closed cover, i.e. $U\subseteq
Y:=\bigcup\limits_{x\in A}V(x),$ and assume without loss of generality that $%
A\subseteq \mathcal{C}(L).$ It follows that $I(Y)=\bigwedge\limits_{x\in
A}x. $ Since $\mathcal{C}(\mathcal{L})$ satisfies the DCC, $%
I(Y)=\bigwedge\limits_{x\in B}x$ for some finite subset $B\subseteq I$. It
follows that%
\begin{equation*}
\overline{Y}\overset{\text{Lemma \ref{Lemma 3.3}}}{=}V(I(Y))=V(\bigwedge%
\limits_{x\in B}x)\overset{\text{\cite[Theorem 2.2]{AbuC}}}{=}%
\bigcup\limits_{x\in B}V(x).
\end{equation*}%
Therefore, $U\subseteq \bigcup\limits_{x\in B}V(x)$ for some finite subset $%
B\subseteq A$.

Conversely, suppose that $x_{1}\geq x_{2}\geq \cdots $ is a descending chain
in $\mathcal{C}(\mathcal{L})$ and consider the induced ascending chain $%
V(x_{1})\subseteq V(x_{2})\subseteq \cdots $. Let $Y=\bigcup\limits_{i=1}^{%
\infty }V(x_{i}).$ By assumption, $Y=\bigcup\limits_{i=1}^{n}V(x_{i})$ for
some $n\in \mathbb{N}$, whence $V(x_{n})=V(x_{m})$ for all $m\geq n$ and
consequently, $x_{n}=x_{m}$ for all $m\geq n$ by Lemma \ref{lemma 1.19}.

\item Assume that $X$ is (countably) compact and suppose that $1=\tilde{%
\bigvee }_{i\in I}x_{i}$ where $x_{i}\in \mathcal{C}(L)$ (and $I$ is
countable). It follows that $\emptyset =V(\tilde{\bigvee }_{i\in
I}x_{i})=\bigcap_{i\in I}V(x_{i})$, i.e. $X=\bigcup_{i\in I}(X\backslash
V(x_{i}))$. Since $X$ is (countably) compact, $X=\bigcup_{j\in
F}(X\backslash V(x_{j}))$ for some finite subset $F$ of $I$ and so $1=\tilde{%
\bigvee }_{j\in F}x_{j}.$ So, $1$ is (countably) finitely constructed. The
converse can be obtained similarly.

\item Let $F\subseteq X$ be a closed irreducible subset. Then $F=V(x)$ for
some $x\in L$, whence $\sqrt{x}\in SI(\mathcal{C}(\mathcal{L}))\subseteq X$
by Proposition \ref{Proposition 3.11}. The uniqueness of the generic point
is obvious.

\item This follows by Lemma \ref{Lemma 3.22}.

\item This follows by applying Corollary \ref{Theorem 3.16} to $(X,\tau
^{cl})=(X,\tau )$.

\item Assume that $X$ is Lindelof (compact). Since $X$ is atomic, $%
X=\bigcup\limits_{p\in Min(X)}V(p)$, whence the open cover $\{V(p)\mid p\in
Min(X)\}$ has a countable (finite) subcover for $X$, i.e. $%
X=\bigcup\limits_{p\in A}V(p)$ for some countable (finite) subset $%
A\subseteq Min(X)$. \textbf{Claim: }$Min(X)=A.$ Let $q\in Min(X).$ Since $%
X=\bigcup\limits_{p\in A}(V(p)$, we have $q\in V(p)$ for some $p\in A$,
whence $q=p$ by the minimality of $q.$ Consequently $Min(X)$ is countable
(finite).

\item This is obtained by applying Proposition \ref{Proposition 3.4} to $%
(X,\tau ^{cl})=(X,\tau )$.
\end{enumerate}
\end{Beweis}

$\blacksquare$

\begin{thm}
\index{$X-top$} \label{Theorem 3.45}Let $X\subseteq L\backslash \{1\}$ and
assume that $\mathcal{L}$ is an $X$-top lattice.

\begin{enumerate}
\item The following are equivalent for the sublattice%
\begin{equation*}
\mathcal{C}^{\prime }(\mathcal{L})=\{x\in \mathcal{C}(L)\mid x%
\tilde{\vee}y=1\text{ and }x\wedge y=\sqrt{0}\text{ for some }y\in \mathcal{C%
}(L)\}
\end{equation*}%
of $\mathcal{C}(\mathcal{L}):$

\begin{enumerate}
\item $(X,\tau )$ is connected.

\item If $x\in L$ is such that $\emptyset \neq V(x)\subsetneq X$, then $V(x)$
is not open.

\item $V(x)\cap V(y)\neq \emptyset $ for all $x\in L$ such that $\sqrt{x}%
\notin \{\sqrt{0},1\}$ and for all $y\in L$ such that $X\backslash
V(x)\subseteq V(y)$.

\item $\mathcal{C}^{\prime }(\mathcal{L})=\{\sqrt{0},1\}$.
\end{enumerate}

\item Let $(X,\tau )$ be $T_{1}$. Then $X$ is singleton if and only if $%
(X,\tau )$ is connected and $\mathcal{C}(\mathcal{L})$ satisfies the
complete max property.

\item If $X$ is coatomic and $Max(X)$ is countable (finite), then $(X,\tau )$
is Lindelof (compact).

\item Let $X$ be coatomic. Then $Max(X)$ is singleton if and only if $%
(X,\tau )$ is connected and each element in $Max(X)$ is completely strongly
irreducible in $(\mathcal{C}(L),\wedge )$ .

\item Let $\mathcal{L}$ be coatomic and $Max(L)\subseteq X$. Then $(X,\tau )$
is ultraconnected if and only if $\mathcal{L}$ is hollow.

\item Let $\emptyset \neq X$ be atomic. Then $(X,\tau )$ is reducible if and
only if $Min(X)=I_{1}\cup I_{2}$ such that $\bigwedge\limits_{p\in
I_{2}}p\nleq q_{l}$ for some $q_{l}\in I_{l}$ and $\bigwedge_{p\in
I_{1}}p\nleq q_{2}$ for some $q_{2}\in I_{2}$ .

\item Let $\emptyset \neq X$ be atomic. Then $(X,\tau )$ is connected if and
only if for every $\emptyset \neq \mathbf{m}\subsetneq Min(X)$\textbf{\ }%
there exists some $q\in X$ such that
\begin{equation*}
\wedge _{p\in \mathbf{m}}p\bigvee \wedge _{p\in Max(X)\backslash \mathbf{m}%
}p\leq q.
\end{equation*}
\end{enumerate}
\end{thm}

\begin{Beweis}
Let $X\subseteq L\backslash \{1\}$ and assume that $\mathcal{L}$ is an $X$%
-top lattice.

\begin{enumerate}
\item Let $x,y\in \mathcal{C}^{\prime }$. Then there are $x^{\prime
},y^{\prime }\in \mathcal{C}(L)$ such that $x\tilde{\vee}x^{\prime }=1$, $%
x\wedge x^{\prime }=\sqrt{0}$, $y\tilde{\vee}y^{\prime }=1$ and $y\wedge
y^{\prime }=\sqrt{0}$. One can check that $x\wedge y$ and $x\tilde{\vee}y$
are also in $\mathcal{C}^{\prime }$ with the corresponding elements $%
x^{\prime }\tilde{\vee}y^{\prime }$ and $x^{\prime }\wedge y^{\prime }$
respectively (recall that if $\mathcal{L}$ is $X$-top then $\mathcal{C}(%
\mathcal{L})$ is distributive by \cite[Theorem 2.2]{AbuC}).

The equivalence $(a)\Leftrightarrow (b)$ is trivial.

$(a\Rightarrow c)$ Let $x,y\in L$ be such that $\sqrt{x}\notin \{\sqrt{0}%
,1\} $ and $X\backslash V(x)\subseteq V(y)$. It follows that $V(x)\cup
V(y)=X $, whence $V(x)\cap V(y)\neq \emptyset $ (otherwise, $X$ will be
disconnected).

$(c\Rightarrow b)$ Suppose that $\emptyset \neq V(x)\subsetneq X$ is open
for some $x\in L$, so that $\sqrt{x}\notin \{\sqrt{0},1\}.$ Let $y\in L$ be
such that $X\backslash V(x)=V(y).$ By our assumption, $V(x)\cap V(y)\neq
\emptyset $ (a contradiction).

$(c\Rightarrow d)$ Let $x\in \mathcal{C}^{\prime }(\mathcal{L})$. Then there
is $y\in \mathcal{C}(L)$ such that $x\wedge y=\sqrt{0}$ and $x\tilde{\vee}%
y=1 $. Clearly, $x$ and $y$ satisfy the conditions stated in $(c),$ whence $%
V(x\tilde{\vee}y)=V(x)\cap V(y)\neq \emptyset $, i.e. $x\tilde{\vee}y\neq 1$%
, which is a contradiction.

$(d\Rightarrow a)$ Suppose that $V(x)\cup V(y)=X,$ $V(x)\cap V(y)=\emptyset $
for some $x,y\in L,$ and assume without loss of generality that $x,y\in
\mathcal{C}(L)$. It is easy to show that $x,y\in \mathcal{C}^{\prime }(%
\mathcal{L})$, and it follows by $(d)$ that $V(x)=X$ or $V(x)=\emptyset $.

\item Let $(X,\tau )$ be $T_{1}.$ If $\mathcal{C}(\mathcal{L})$ satisfies
the complete max property, then applying Corollary \ref{Theorem 3.16} to $(X,\tau
)=(X,\tau ^{cl}),$ we conclude that $X$ is discrete. If $X$ is moreover
connected, then $X$ is indeed a singleton. The converse is trivial.

\item Let $X$ be coatomic and $Max(X)$ be countable (finite). Let $\mathcal{A%
}=\{X\backslash V(x)\mid x\in A\}$ be an open cover for $X$. Then $%
\bigcap\limits_{x\in A}V(x)=\emptyset $ and so for any $p\in Max(X)$, there
exists $x_{p}\in A$ such that $p\notin V(x_{p}).$ \textbf{Claim: }$%
\bigcap\limits_{p\in Max(X)}V(x_{p})=\emptyset $. Suppose that $q\in
\bigcap\limits_{p\in Max(X)}V(x_{p})$. Since $X$ is coatomic, $q\leq \tilde{p%
}$ for some $\tilde{p}\in Max(X)$ and so $\tilde{p}\in \bigcap\limits_{p\in
Max(X)}V(x_{p})$, a contradiction.

It follows that $X=\bigcup\limits_{p\in Max(X)}(X\backslash V(x_{p})),$ i.e.
$\{X\backslash V(x_{p}\mid p\in Max(X)\}$ is a countable (finite) subcover
of $\mathcal{A}$ for $X.$

\item Let $X$ be coatomic.

$(\Rightarrow )$ Assume that $Max(X)=\{p\}$. For all $q\in X$, $q\leq p$ as $%
X$ is coatomic and so $p$ is completely irreducible in the $(\mathcal{C}%
,\wedge )$. Also, if $X=V(x)\cup V(y)$ and $V(x),V(y)\neq \emptyset $, then $%
p\in V(x)\cap V(y)$ and so $X$ is connected.

$(\Leftarrow )$ Suppose that $|Max(X)|\geq 2$ and let $Max(X)=\mathbf{%
M^{\prime }}\cup \mathbf{M^{\prime \prime }}$ with $\mathbf{M^{\prime }}\cap
\mathbf{M^{\prime \prime }}=\emptyset $ for some $\emptyset \neq \mathbf{%
M^{\prime }}\subsetneq Max(X)$. Set%
\begin{equation*}
A:=\{p\in X\mid p\leq q\text{ for some }q\in \mathbf{M^{\prime }}\text{ and }%
p\nleq q\text{ }\forall \text{ }q\in \mathbf{M^{\prime \prime }}\},
\end{equation*}%
$B:=X\backslash A,$ $x:=\bigwedge\limits_{p\in A}p$ and $y:=\bigwedge%
\limits_{p\in B}p$.

\textbf{Claim:} $V(x)\cap V(y)=\emptyset $.

Suppose that $\tilde{p}\in V(x)\cap V(y)$, whence $y\leq \tilde{p}\leq
\tilde{q}$ for some $\tilde{q}\in Max(X)\mathbf{.}$ Since $\tilde{q}$ is
completely strongly irreducible, $\tilde{q}\in \mathbf{M^{\prime }}$:
otherwise, $\tilde{q}\in \mathbf{M^{\prime \prime }}$ and $%
x=\bigwedge\limits_{p\in A}p\leq \tilde{q}$ implies that $p^{\prime }\leq
\tilde{q}\in \mathbf{M^{\prime \prime }}$ for some $p^{\prime }\in A$, a
contradiction. Hence, $y\leq \tilde{q}\in \mathbf{M^{\prime }}$. Similarly,
since $\tilde{q}$ is completely strongly irreducible, $q^{\prime }\leq
\tilde{q}$ for some $q^{\prime }\in B$, which is a contradiction. Therefore $%
V(x)\cap V(y)=\emptyset $, and $V(x)$ and $V(y)$ are non-empty ($\mathbf{%
M^{\prime }}\subseteq V(x)$ and $\mathbf{M^{\prime \prime }}\subseteq V(y)$)
with $V(x)\cup V(y)=X$, whence $X$ is disconnected.

\item Let $\mathcal{L}$ be coatomic and $Max(L)\subseteq X.$

$(\Rightarrow )$ Assume that $X$ is ultraconnected. Let $x,y\in L\backslash
\{1\}$. Since $\mathcal{L}$ is coatomic, there are $p,q\in Max(\mathcal{L}%
)\subseteq X$ with $x\leq p$ and $y\leq q,$ whence $V(x)$ and $V(y)$ are
non-empty. By assumption, $X$ is ultraconnected, whence $V(x\vee y)=V(x)\cap
V(y)\neq \emptyset .$ Hence $x\vee y\neq 1.$ Consequently, $\mathcal{L}$ is
hollow.

$(\Leftarrow )$ Assume that $\mathcal{L}$ is hollow. Let $V(x)$ and $V(y)$
be non-empty closed subsets for some $x,y\in L$. Then $x,y\in L\backslash
\{1\}$, whence $x\vee y\neq 1$ as $\mathcal{L}$ is hollow. Since $L$ is
coatomic, $x\vee y\leq q$ for some $q\in Max(\mathcal{L})\subseteq X$. Hence
$V(x)\cap V(y)=V(x\vee y)\neq \emptyset $. Therefore, $X$ is ultraconnected.

\item Let $X$ be reducible, i.e. $X=V(x)\cup V(y)$ for some $x,y\in L$ such
that $V(x)\subsetneq X$ and $V(y)\subsetneq X$. Set
\begin{equation*}
I_{1}=\{p\in Min(X)\mid x\leq p\}\text{ and }I_{2}=\{p\in Min(X)\mid y\leq
p\}.
\end{equation*}%
Since $X$ is atomic, $\sqrt{x}=\bigwedge\limits_{p\in I_{1}}p$ and $\sqrt{y}%
=\bigwedge\limits_{p\in I_{2}}p$. Indeed, $\sqrt{x}\nleq q_{2}$ for some $%
q_{2}\in I_{2}:$ otherwise, $\sqrt{x}\leq p$ $\forall $ $p\in I_{2}$ and it
follows that $V(x)=X$. Similarly, $\sqrt{y}\nleq q_{1}$ for some $q_{1}\in
I_{1}$. The converse is trivial.

\item Let $\emptyset \neq X$ be atomic.

$(\Rightarrow )$ Assume that $X$ is connected. Let $\emptyset \neq \mathbf{m}%
\subsetneq Min(X),$ $x:=\wedge _{p\in \mathbf{m}}p$ and $y=\wedge _{p\in
Max(X)\backslash \mathbf{m}}p$. Since $X$ is atomic, $X=V(x)\cup V(y)$.
Since $X$ is connected, $V(x\vee y)=V(x)\cap V(y)\neq \emptyset ,$ i.e. $%
\exists $ $q\in X$ such that $x\vee y\leq q$.

$(\Leftarrow )$ Suppose that $X=V(x)\cup V(y)$ for some $x,y\in L.$ Set
\begin{equation*}
\mathbf{m}^{\prime }:=\{p\in Min(X)\cap V(x)\}\text{ and }\mathbf{m^{\prime
\prime }}:=Min(X)\backslash \mathbf{m^{\prime }}.
\end{equation*}

\textbf{Case 1: }$\mathbf{m^{\prime }}=\emptyset $. In this case, $X=V(y)$.

\textbf{Case 2:} $\mathbf{m^{\prime }}=Min(X).$ In this case, $X=V(x)$.

\textbf{Case 3:} $\emptyset \neq \mathbf{m^{\prime }}\subsetneq Min(X).$ By
our assumption, $\sqrt{x}\vee \sqrt{y}\leq q$ for some $q\in X$ and so
\begin{equation*}
V(x)\cap V(y)=V(\sqrt{x})\cap V(\sqrt{y})=V(\sqrt{x}\vee \sqrt{y})\neq
\emptyset .
\end{equation*}%
Consequently, $X$ is connected.
\end{enumerate}
$\blacksquare$
\end{Beweis}

\begin{ex}
Let $M$ be a left module over an arbitrary ring $R$. Consider $X_1 =
Spec^p(M)$ and $X_2 = Spec^c(M)$. Suppose that $\sqrt{0}=0$ (e.g. the $%
\mathds{Z}$-module $\mathds{Z}[x]$). Then the set $\mathcal{C}^{\prime}$
which was described in Theorem \ref{Theorem 3.45} (1) is the set of the
prime radical direct summands (resp. the coprime radical direct summands).
\end{ex}

\begin{cor}
\label{Corollary 3.46}Let $X\subseteq L\backslash \{1\}$ and assume that $%
\mathcal{L}$ is an $X$-top lattice.

\begin{enumerate}
\item Let $X$ be atomic, coatomic with $|Max(X)|\leq |Min(X)|$ and $V(p)$ is
open $\forall $\textbf{\ }$p\in Min(X),$ then $(X,\tau )$ is Lindelof
(compact) if and only if $Max(X)$ is countable (finite).

\item Let $X=Max(\mathcal{L})$. Then $\left\vert Max(\mathcal{L})\right\vert
=1$ if and only if $(X,\tau )$ is connected and $\mathcal{C}(\mathcal{L})$
satisfies the complete max property.
\end{enumerate}
\end{cor}

\begin{Beweis}
\begin{enumerate}
\item If $(X,\tau )$ is Lindelof, then $Min(X)$ is countable by Remark \ref%
{Remark 3.44}(10). Conversely, assume that $Max(X)$ is countable (finite).
Let $\mathcal{O}=\{X\backslash V(x)\mid x\in A\subseteq L\}$ be an open
cover for $X$, i.e. $\bigcap\limits_{x\in A}V(x)=\emptyset $ and assume
without loss of generality that $V(x)\neq \emptyset $ for each $x\in A$ (If $%
V(y)=\emptyset $ for some $y\in A,$ then $\{X\backslash V(y)\}$ is a finite
subcover of $X$). Pick $x^{\prime }\in A$ and set $\mathbf{M}:=\{q\in
Max(X)\mid x^{\prime }\leq q\}$. Observe that $\mathbf{M}$ is non-empty as $%
V(x^{\prime })\neq \emptyset $ and $X$ is coatomic. For each $q\in \mathbf{M,%
}$ pick $X\backslash V(x_{q})\in \mathcal{O}$ that contains $q$.

\textbf{Claim:} $x^{\prime }\vee \bigvee\limits_{q\in \mathbf{M}}x_{q}\nleq
p $ for each $p\in Max(X).$

\textbf{Case (1): }$p\in \mathbf{M}$. In this case, $x_{p}\nleq p$ and so $%
x^{\prime }\vee \bigvee\limits_{q\in \mathbf{M}}x_{q}\nleq p$.

\textbf{Case (2): }$p\in Max(X)\backslash \mathbf{M}$. In this case, $%
x^{\prime }\nleq p$ and so $x^{\prime }\vee \bigvee\limits_{q\in \mathbf{M}%
}x_{q}\nleq p$.

Therefore, $V(x^{\prime }\vee \bigvee\limits_{q\in \mathbf{M}%
}x_{q})=\emptyset $ and
\begin{equation*}
\{X\backslash V(x^{\prime })\}\cup \{X\backslash V(x_{q})\mid q\in \mathbf{M}%
\}
\end{equation*}%
is a countable (finite) subcover of $\mathcal{O}$ as $Max(X)$ is countable
(finite).

\item Assume that $Max(\mathcal{L})=X$, whence $Max(X)=X=Max(\mathcal{L})$.
It follows by Corollary \ref{Theorem 3.16} that $X$ is $T_{1}$. So, we can use
now Theorem \ref{Theorem 3.45} (2).
\end{enumerate}
\end{Beweis}

$\blacksquare$

\begin{thm}
\index{$X-top$} \label{Proposition 3.48}Let $X\subseteq L\backslash \{1\}$
and assume the $\mathcal{L}$ is an $X$-top lattice.

\begin{enumerate}
\item There is a 1-1 correspondence
\begin{equation*}
\mathcal{C}(L)%
\hspace{0.5cm}\longleftrightarrow \hspace{0.5cm}\text{closed sets in }%
(X,\tau ).
\end{equation*}

\item If $SI(\mathcal{C}(\mathcal{L}))\subseteq X$, then there is a 1-1
correspondence
\begin{equation*}
X\hspace{0.5cm}\longleftrightarrow \hspace{0.5cm}\text{Irreducible closed
sets in }(X,\tau ).
\end{equation*}

\item If $SI(\mathcal{C}(\mathcal{L}))\subseteq X$, then there is a 1-1
correspondence
\begin{equation*}
Min(X)\hspace{0.5cm}\longleftrightarrow \hspace{0.5cm}\text{Irreducible
components in }(X,\tau ).
\end{equation*}
\end{enumerate}
\end{thm}

\begin{Beweis}
Since $\mathcal{L}$ is $X$-top, the set of closed sets in $X$ is given by $%
\mathcal{V}=\{V(y)\mid y\in L\}.$ Define%
\begin{eqnarray*}
f &:&\mathcal{C}(L)\longrightarrow \mathcal{V},\text{ }x\mapsto V(x); \\
g &:&\mathcal{V}\longrightarrow \mathcal{C}(L),\text{ }V(y)\mapsto \sqrt{y}.
\end{eqnarray*}

\begin{enumerate}
\item For any $x\in \mathcal{C}(L)$ and $y\in L,$ we have%
\begin{eqnarray*}
(g\circ f)(x) &=&g(V(x))=\sqrt{x}=x; \\
(f\circ g)(V(y)) &=&f(\sqrt{y})=V(\sqrt{y})=V(y).
\end{eqnarray*}%
So, $f$ provides a 1-1 correspondence $\mathcal{C}(L)\longleftrightarrow
\mathcal{V}$ with inverse $g.$

\item Consider the restrictions of $f$ to $X$ and of $g$ to the class of
irreducible closed varieties. For every $x\in X$, the variety $V(x)$ is irreducible by Proposition \ref%
{Proposition 3.4} (2). On the other hand, if $V(y)$ is irreducible for some $%
y\in L$, then $\sqrt{y}$ is strongly irreducible in $\mathcal{C}(\mathcal{L}%
) $ by Proposition \ref{Proposition 3.11}, whence $\sqrt{y}\in X$ by our
assumption.

\item Consider the restrictions of $f$ to $Min(X)$ and of $g$ to the class
of irreducible components in $(X,\tau ).$

For every $x\in Min(X).$ By (2), $V(x)$ is irreducible. Suppose that $%
V(x)\subseteq V(y)$ for some $y\in L$ with $V(y)$ irreducible. Since $SI(%
\mathcal{C}(\mathcal{L}))\subseteq X,$ it follows by (2) that $\sqrt{y}\in X$%
, whence $\sqrt{y}\leq x.$ However, $x$ in minimal in $X,$ whence $x=\sqrt{y}
$ and $V(x)=V(y).$

On the other hand, let $A$ be an irreducible component in $(X,\tau ).$ Any
irreducible component is closed. Moreover, $I(A)$ is strongly irreducible in
$\mathcal{C}(\mathcal{L})$ as $A$ is irreducible, hence $I(A)\in X$. Suppose
that $p\leq I(A)$ for some $p\in X.$ It follows that $A=\overline{A}%
=V(I(A))\subseteq V(p)=\overline{\{p\}}$. However, $V(p)$ is irreducible as
it is the closure of a singleton, so $V(p)=A$ as $A$ is an irreducible
component. So, $p=I(A)$. Consequently, $I(A)\in Min(X).$
\end{enumerate}
\end{Beweis}

$\blacksquare$

\begin{ex}
The first correspondence ($\mathcal{C}(\mathcal{L}(M))\hspace{0.5cm}%
\longleftrightarrow \hspace{0.5cm}$closed sets in $(X,\tau )$) of Theorem %
\ref{Proposition 3.48} holds for any $X\subseteq \mathcal{L}(M)\backslash
\{M\}$ such that $\mathcal{L}$ is $X$-top, as well as for any $X\subseteq
\mathcal{L}(M)\backslash \{0\}$ such that $\mathcal{L}$ is dual $X$-top. So,
this result recovers \cite[4.12 and 5.16]{Abu}, \cite[3.27]{Abu2011-a} and
\cite[3.23]{Abu2011-b} as special cases.
\end{ex}

The following table summarizes some of the results we obtained in this
section. Some of them generalize results in the literature on Zariski-like
topologies for left modules over associative rings, which can be recovered
now as special cases. At several occasions, our results were obtained under
conditions and assumptions weaker than those in the corresponding results in
the literature.
\begin{table}[tbp]
\centering
\par
\begin{adjustbox}{width=1\textwidth,height=0.48\textheight}
\begin{tabular}{ | m{8em} | m{8em} | m{8em}| m{8em} | }
 \hline
 \textbf{Assumption $\&$ location} & $X$\textbf{-top lattice} $\mathcal{L}$ & $(X,\tau)$ & \textbf{Results recovered} \\ \hline
   Proposition \ref{Proposition 3.15}  & $Max(X) = X$ & $T_1$ & \cite[4.27, 5.33]{Abu}, \cite[3.45]{Abu2011-a})\\ \hline
 Proposition \ref{Proposition 3.24}    & $\mathcal{C}(\mathcal{L})$ satisfies the ACC & Each set in $X$ is compact  & \\ \hline
   Remark \ref{Remark 3.44} (3)  & $\mathcal{C}(\mathcal{L})$ satisfies the ACC &  Noetherian &  \cite[4.12, 5.16]{Abu}) \\ \hline
    Proposition \ref{Proposition 3.24} & $\mathcal{C}(\mathcal{L})$ satisfies the ACC & Each open set in $X$ is compact  & \\ \hline
   Remark \ref{Remark 3.44} (4) & $\mathcal{C}(\mathcal{L})$ satisfies the DCC & Artinian & \cite[4.12, 5.16]{Abu}\\ \hline
    Remark \ref{Remark 3.44} (4)  & $\mathcal{C}(\mathcal{L})$ satisfies the DCC & Every closed cover for any subset of $X$ has a finite subcover & \\ \hline
     Theorem \ref{Theorem 3.16 000} &  $Max(X) = X$ and $\mathcal{C}(\mathcal{L})$ satisfies the complete max property& Discrete & \cite[4.28, 5.34]{Abu}, \cite[3.46]{Abu2011-a}, \cite[4.33]{AbuC})\\ \hline
     Theorem \ref{Theorem 3.45} (1) & $\mathcal{C}' = \{\sqrt{0}, 1\}$ & Connected & \\ \hline

 Corollary \ref{Corollary 3.12} & $I(A)\in SI(\mathcal{C}(\mathcal{L}))$ & $A\subseteq X$ is irreducible & \cite[3.30, 3.31]{Abu2011-a}\\ \hline
Corollary \ref{Corollary 3.12}& $I(A)$ is irreducible in $\mathcal{C}(\mathcal{L})$ & $A\subseteq X$ is irreducible & \cite[3.30, 3.31]{Abu2011-a} \\ \hline
Corollary \ref{Corollary 3.12}& $\sqrt{0}$ is irreducible in $\mathcal{C}(\mathcal{L})$ & irreducible & \cite[3.30, 3.31]{Abu2011-a} \\ \hline
   $SI(\mathcal{C}(\mathcal{L}))\subseteq X$ (\ref{Proposition 3.48}) & $\sqrt{x}\in X$ & $V(x)$ is irreducible & \cite[4.17, 5.22]{Abu}, \cite[3.27]{Abu}, \cite[3.33]{Abu2011-a}, \cite[3.6]{BH2008-a}, \cite[4.28]{AbuC}\\ \hline
    $SI(\mathcal{C}(\mathcal{L}))\subseteq X$ (\ref{Proposition 3.48}) & $\sqrt{x}\in Min(X)$ & $V(x)$ is irreducible component & \cite[5.22]{Abu}, \cite[4.17]{Abu}, \cite[3.27]{Abu2011-a}, \cite[3.33]{Abu2011-a}, \cite[4.28]{AbuC}\\ \hline
      $Max(\mathcal{L}) = X$ and $\mathcal{C}(\mathcal{L})$ satisfies the complete max property (\ref{Theorem 3.45} (2))& $|Max(X)|=1$ & Connected & \\ \hline
    Remark \ref{Remark 3.44}  (5)  & $1$ is finitely constructed  & Compact & \\ \hline
   Remark \ref{Remark 3.44} (5)   & $1$ is countably  constructed  & Lindelof & \\ \hline

  \end{tabular}
\end{adjustbox}
\caption{Examples on the Interplay between topological properties of $(X,%
\protect\tau )$ and algebraic properties of the $X-top$ lattice $\mathcal{L}$%
}
\label{Interplay}
\end{table}

\newpage

\begin{lem}
\index{$LAT(M)$} \label{Remark 3.47}Let $R$ be a ring and $M$ a top$^{p}$%
-module, that is $\mathcal{L}:=LAT(_{R}M)$ is $Spec^{p}(M)$-top. Then
\begin{equation*}
SI(\mathcal{C}(LAT(M)))\subseteq Spec^{p}(M).
\end{equation*}
\end{lem}

\begin{Beweis}
Let $N$ be strongly irreducible in $\mathcal{C}(\mathcal{L}).$ Suppose that $%
IK\subseteq N$ for some ideal $I\leq R$ and a submodule $K\leq M$. Then $%
IK\subseteq P$ for any prime submodule $P\in V(N)$, whence $IM\subseteq P$
or $K\subseteq P$ and so $%
\sqrt{IM}\subseteq P$ or $\sqrt{K}\subseteq P$, whence $\sqrt{IM}\cap \sqrt{K%
}\subseteq P$ for all $P\in V(N)$. Since $N$ is radical, $\sqrt{IM}\cap
\sqrt{K}\subseteq N.$ By assumption, $N$ is strongly irreducible in $%
\mathcal{C}(\mathcal{L})$, whence $IM\subseteq \sqrt{IM}\subseteq N$ or $%
K\subseteq \sqrt{K}\subseteq N$. Therefore, $N\in Spec^{p}(M).$
\end{Beweis}

$\blacksquare$

\begin{ex}
\label{ex-p}Let $R$ be a ring and $M$ a top$^{p}$-module. By Lemma \ref%
{Remark 3.47}, we have $SI(\mathcal{C}(LAT(_{R}M)))\subseteq Spec^{p}(M).$
So, all the 1-1 correspondences in Theorem \ref{Proposition 3.48} hold for
this special case. Behboodi and Haddadi proved the second correspondence in
\cite[Corollary 3.6]{BH2008-a}.
\end{ex}

\begin{ex}
Let $R$ be a ring and $_{R}M$ a left top$^{c}$-module (i.e. $\mathcal{L}%
=LAT(_{R}M)$ is $X$-top, where $X=Spec^{c}(M)$). If $_{R}M$ is completely
distributive, then $SI(\mathcal{C}(\mathcal{L}))\subseteq X$ by \cite[%
Proposition 5.19]{Abu} and the 1-1 correspondences of Theorem \ref%
{Proposition 3.48} hold. In \cite[Proposition 5.22]{Abu}, these
correspondences were proved under the additional condition that every
coprime submodule of $M$ is strongly irreducible.
\end{ex}

\begin{ex}
Let $R$ be a ring and $_{R}M$ a left top$^{s}$-module (\emph{i.e.} $\mathcal{%
L}=LAT(_{R}M)$ is dual $X$-top, where $X=Spec^{s}(M)$). By \cite[Proposition
4.14]{Abu}, $SH(\mathcal{H}(\mathcal{L}))\subseteq X$ and so the 1-1
correspondences of Theorem \ref{Proposition 3.48} hold. These were proved in
this special case in \cite[Proposition 4.17]{Abu} under the additional
condition that every second submodule of $M$ is strongly hollow.
\end{ex}

\begin{ex}
Let $R$ be a ring and $_{R}M$ a left top$^{fp}$-module (\emph{i.e.} $%
\mathcal{L}=LAT(_{R}M)$ is $X$-top, where $X=Spec^{fp}(M)$). If $_{R}M$ is
duo, then $SI(\mathcal{C}(\mathcal{L}))\subseteq X$ by \cite[3.30]{Abu2011-a}
and the 1-1 correspondences of Theorem \ref{Proposition 3.48} hold. These
were also obtained under the same condition in \cite[Proposition 3.33]%
{Abu2011-a}.
\end{ex}

\begin{ex}
Let $R$ be a ring and $_{R}M$ a left top$^{fc}$-module (\emph{i.e.} $%
\mathcal{L}=LAT(_{R}M)$ is dual $X$-top, where $X=Spec^{fc}(M)$). If $_{R}M$
is duo, then $SH(\mathcal{H}(L))\subseteq X$ by \cite[Proposition 3.25]%
{Abu2011-b} and Proposition \ref{Proposition 3.11} and the 1-1
correspondences of Theorem \ref{Proposition 3.48} hold. These were also
obtained under the same condition for this special case in \cite[Proposition
3.28]{Abu2011-b}.
\end{ex}

\begin{ex}
\label{ex-f}Let $R$ be a ring and $_{R}M$ a left top$^{f}$-module (\emph{i.e.%
} $\mathcal{L}=LAT(_{R}M)$ is dual $X$-top, where $X=Spec^{f}(M)$). If $%
_{R}M $ has the property that $H(A)$ is first whenever $A$ is irreducible,
then $SH(\mathcal{H}(L))\subseteq X$ and so the 1-1 correspondences of
Theorem \ref{Proposition 3.48} hold. This was proved under the same
condition in \cite{AbuC}.
\end{ex}

\begin{ex}
\label{Example 3.49}Let $R$ be a PID with an infinite number of non-zero
prime ideals (e.g. $R=\mathbb{Z}$), $\mathcal{L}:=Ideal(R),$ $X=Max(R)$ and
consider the topological space $(X,\tau ).$

\begin{enumerate}
\item $X=V(0)$ is irreducible since $0$ is strongly irreducible. However, $0=%
\sqrt{0}\notin X$ and so $X$ is not sober by Remark \ref{Remark 3.44} (7),
whence not spectral.

\item $X$ is $T_{1}$ as $Max(X)=X$.

\item $X$ is cofinite: consider a closed set $\emptyset \neq V(I)\subsetneq
X,$ where $I=(a)$ for some $a\in R\backslash \{0\}$. Since $R$ is a PID, the
unique prime factorization of $a$ implies that $I$ is contained in a finite
number of primes, i.e. $V(I)$ is finite.

\item $X$ is not regular, not $T_{2},$ and not normal. Observe that $X$ is
infinite and cofinite, so it does not have disjoint non-empty open sets,
although it has disjoint non-empty closed sets.
\end{enumerate}
\end{ex}

\begin{ex}
Let $R$ be a ring, $M$ a left $R$-module, $X\subseteq LAT(_{R}M)\backslash \{M\}$
(resp. $X\subseteq LAT(_{R}M)\backslash \{0\}$) and assume that $\mathcal{L}:=LAT(_{R}M)$
is $X$-top (resp. dual $X$-top). If $\mathcal{C}(\mathcal{L})$ is uniform (resp. $\mathcal{H}(\mathcal{L})$
is hollow), then $(X,\tau )$ (resp. $(X,\tau ^{0}$)) is connected by Theorem %
\ref{Theorem 3.45} (1).
\end{ex}

\begin{ex}
Let $R$ be a commutative domain, $\mathcal{L}:=Ideal(R),$ $X\subseteq
Ideal(R)\backslash \{R\}$ (resp. $X\subseteq Ideal(R)\backslash \{0\}$), and
assume that $\mathcal{L}$ is $X$-top (resp. $\mathcal{L}$ is dual $X$-top).
If $\sqrt{0}=0$ (resp. $\sum\limits_{p\in X}p=R$), Then $(X,\tau )$ (resp. $%
(X,\tau ^{0})$) is connected.
\end{ex}

\begin{ex}
Let $R$ be a UFR with zero devisors. Consider $\mathcal{L}:=Ideal(R),$ $%
X:=Spec(R)$ (the prime spectrum of $R$) and assume that $Min(X)$ is infinite
(e.g. $R=Z_{n}[x]$ with $n$ not prime). Notice that $\sqrt{0}=0$ (since $%
Min(X)$ is infinite, if $0\neq x\in \sqrt{0}$ then $x\in
\bigcap\limits_{p\in Min(X)}q,$ but this is impossible as $R$ is a UFR).

\begin{itemize}
\item $(X,\tau )$ is connected by Theorem \ref{Theorem 3.45} (7).

\textbf{Claim: }the intersection of any infinite collection of minimal
elements of $X$ is zero. Suppose that $0\neq I:=\bigcap\limits_{q\in \mathbf{%
m}^{\prime }}q$ for some infinite subcollection $\mathbf{m}^{\prime }$ of $%
Min(R)$. For any $x\in I\backslash \{0\}$, we have $x=p_{1}\cdots p_{n}$
where $p_{1},\cdots ,p_{n}$ are prime elements of $R.$ Notice that $%
p_{1},\cdots ,p_{n}\in I.$ For every $q\in \mathbf{m}^{\prime },$ we have $%
q=(p_{i})$ for some $i\in \{1,2,\cdots ,n\}$, whence $\mathbf{m}^{\prime }$
is finite (a contradiction).

\item $(X,\tau )$ is reducible by Remark \ref{Remark 3.44} (1). To prove
this, suppose that $(X,\tau )$ is irreducible and that $I\cap J=0$ for some
ideals $I,J\leq R$. Then $V(0)=V(I\cap J)=V(I)\cup V(J)=V(\sqrt{I})\cup V(%
\sqrt{J})=V(\sqrt{I}\cap \sqrt{J})$, whence $\sqrt{I}\cap \sqrt{J}=\sqrt{0}%
=0.$ Since $\sqrt{0}\in SI(\mathcal{C}(\mathcal{L}))$ (by Remark \ref{Remark
3.44} (1)), it follows that $I=\sqrt{I}=0$ or $J=\sqrt{J}=0$, whence $R$ is
a domain, a contradiction.
\end{itemize}
\end{ex}

\begin{ex}
\label{Example 3.51}Let $(G,+)$ be a group and set \vspace{15pt}

\hspace{15pt} $L :=\{H\mid H\trianglelefteq G\text{ is a normal subgroup of }%
G\}$,

\hspace{15pt} $X :=\{H\mid H\trianglelefteq G\text{ is a finite normal
subgroup of }G\}\backslash \{G\}$. \vspace{15pt}

Notice that $\mathcal{L}=(L,\cap ,+,G,0)$ is a complete lattice endowed with
$\bigvee\limits_{i\in I}N_{i}:=\sum\limits_{i\in I}N_{i}$ and $%
\bigwedge\limits_{i\in I}N_{i}:=\bigcap\limits_{i\in I}N_{i}.$

\begin{enumerate}
\item $\mathcal{C}(L)=X\cup \{G\}$ as the intersection of any non-empty
family of finite normal subgroups is a finite normal subgroup.

\item $SI(\mathcal{C}(\mathcal{L}))\subseteq X$ and so all the 1-1
correspondences of Theorem \ref{Proposition 3.48} hold.

\item $0=\sqrt{0}\in X$ and so $(X,\tau ^{cl})$ is irreducible and connected
(observe that $\overline{\{0\}}=X$ and $\overline{\{0\}}$ is irreducible).

\item $\mathcal{C}(\mathcal{L})$ satisfies the DCC but need not satisfy the
ACC (e.g. a $p$-quasicyclic group \cite{JG2013}).

\item $SI(\mathcal{C}(\mathcal{L}))=X$ if and only if $\mathcal{L}$ is an $X$%
-top lattice.

\item If $\mathcal{L}$ is $X$-top, then the intersection of any nonzero
elements in $X$ is nonzero.

\item By Corollary \ref{Theorem 3.16}: $(X,\tau ^{cl})$ is $T_{1}$ $%
\Leftrightarrow $ $(X,\tau ^{cl})$ is a singleton $\Leftrightarrow $ $%
(X,\tau ^{cl})$ is $T_{2}$ $\Leftrightarrow $ $(X,\tau ^{cl})$ is discrete.

\item Suppose that $\mathcal{L}$ is an $X$-top lattice and $(X,\tau ^{cl})$
is compact with each element in $G$ having a finite order. Then $G$ is a
finite $p$-group for some prime $p$. Indeed, since $X$ is compact, by
Theorem \ref{Theorem 3.45} (5), $G$ is finitely constructed. But $G$ is the
union of all proper cyclic subgroups, say $G=\sum\limits_{i\in I}H_{i}$.
Then $G=\sum\limits_{j\in F}H_{j}$ where $F$ is a finite subset of $I$.
Hence $G$ is finite. Consequently, the Pr\"{u}fer group is not $X$-top ($X$ is the set of all
proper subgroups) as it is infinite.
\end{enumerate}
\end{ex}

\begin{ex}
\label{Example 3.53}Let $(G,+)$ be a group, $Z(G)$ the center of $G$ and set%
\vspace{15pt}

\hspace{15pt} $L :=\{H\mid H\trianglelefteq G\text{ is a normal subgroup of }%
G\}$,

\hspace{15pt} $X :=\{H\mid H\leq Z(G)\}\backslash \{G\}.$ \vspace{15pt}

Notice that $\mathcal{L}=(L,\cap ,+,G,0)$ is a complete lattice with $%
\bigvee\limits_{i\in I}N_{i}:=\sum\limits_{i\in I}N_{i}$ and $%
\bigwedge\limits_{i\in I}N_{i}:=\bigcap\limits_{i\in I}N_{i}.$

\begin{enumerate}
\item $\mathcal{C}(L)=X\cup \{G\}$ as the intersection of any non-empty
family of subgroups of the center is again in the center.

\item $SI(\mathcal{C}(\mathcal{L}))\subseteq X$ and so all correspondences
of Theorem \ref{Proposition 3.48} hold.

\item $0=\sqrt{0}\in X$ and so $(X,\tau ^{cl})$ is irreducible and connected.

\item By Corollary \ref{Theorem 3.16}: $(X,\tau ^{cl})$ is $T_{1}$ $%
\Leftrightarrow $ $(X,\tau ^{cl})$ is singleton $\Leftrightarrow $ $(X,\tau
^{cl})$ is $T_{2}$ $\Leftrightarrow $ $(X,\tau ^{cl})$ is discrete.

\item $SI(\mathcal{C}(\mathcal{L}))=X$ $\Leftrightarrow $ $\mathcal{L}$ is $%
X $-top. Hence, if $\mathcal{L}$ is an $X$-top lattice, then the
intersection of any distinct nonzero subgroups in $X$ is nonzero.

\item If $G$ is finite, then $(X,\tau ^{cl})$ is spectral by Remark \ref%
{Remark 3.21}.

\item Suppose that $\mathcal{L}$ is an $X$-top lattice and $(X,\tau ^{cl})$
is compact with each element in $G$ having a finite order. Then $G$ is a
finite $p$-group for some prime $p$.

\item $X$ is coatomic and $Z(G)$ is the unique maximal element of $X$.

\item If $\mathcal{L}$ is $X$-top, then $X$ is compact as $X$ is coatomic
and $Max(X)$ is finite (by Theorem \ref{Theorem 3.45} (3)).
\end{enumerate}
\end{ex}

\begin{ex}
Let $(G,+)$ be a group, $Z(G)$ the center of $G$ and set\vspace{15pt}

\hspace{15pt} $L :=\{H\mid H\trianglelefteq G\text{ is a normal subgroup of }%
G\}$,

\hspace{15pt} $X :=\{H\mid H\leq Z(G) \text{ is finite }\}\backslash \{G\}.$
\vspace{15pt}

Notice that $\mathcal{L}=(L,\cap ,+,G,0)$ is a complete lattice with $%
\bigvee\limits_{i\in I}N_{i}:=\sum\limits_{i\in I}N_{i}$ and $%
\bigwedge\limits_{i\in I}N_{i}:=\bigcap\limits_{i\in I}N_{i}.$

\begin{enumerate}
\item $\mathcal{C}(L)=X\cup \{G\}$ as the intersection of any non-empty
family of finite subgroups of the center is again finite and in the center.

\item $SI(\mathcal{C}(\mathcal{L}))\subseteq X$ and so all correspondences
of Theorem \ref{Proposition 3.48} hold.

\item $0=\sqrt{0}\in X$ and so $(X,\tau ^{cl})$ is irreducible and connected.

\item By Corollary \ref{Theorem 3.16}: $(X,\tau ^{cl})$ is $T_{1}$ $%
\Leftrightarrow $ $(X,\tau ^{cl})$ is singleton $\Leftrightarrow $ $(X,\tau
^{cl})$ is $T_{2}$ $\Leftrightarrow $ $(X,\tau ^{cl})$ is discrete.

\item $SI(\mathcal{C}(\mathcal{L}))=X$ $\Leftrightarrow $ $\mathcal{L}$ is $%
X $-top. Hence, if $\mathcal{L}$ is an $X$-top lattice, then the
intersection of any distinct nonzero subgroups in $X$ is nonzero and so $X$
can only be $\{0\}$ or a collection of $p$-subgroups for some fixed prime $p$%
. Otherwise, $H\in X$ has order $p^{n}q^{m}l$ with primes $p$ and $q$ not
dividing $l$ and so by the Sylow Theorem \cite[Theorem 5.2]{JG2013} there is
a Sylow $p$-subgroup $K_{1}$ of order $p^{n}$ and a Sylow $q$-subgroup $%
K_{2} $ of order $q^{m}$. By Lagrange's Theorem \cite[Theorem 1.26]{JG2013},
the order of their intersection must divide $p^{n}$ and $q^{m}$ and so the
intersection must be zero, whereas $K_{1}$ and $K_{2}$ are nonzero elements
of $X$. The uniqueness of $p$ is clear also by Lagrange's Theorem.

\item If $G$ is finite, then $(X,\tau ^{cl})$ is spectral by Remark \ref%
{Remark 3.21}.

\item Suppose that $\mathcal{L}$ is an $X$-top lattice and $(X,\tau ^{cl})$
is compact with each element in $G$ having a finite order. Then $G$ is a
finite $p$-group for some prime $p$.

\item $X$ is coatomic and $Z(G)$ is the unique maximal element of $X$.

\item If $\mathcal{L}$ is $X$-top, then $X$ is compact as $X$ is coatomic
and $Max(X)$ is finite (by Theorem \ref{Theorem 3.45} (3)).
\end{enumerate}
\end{ex}

\end{document}